# SAMPLE-PATH LARGE DEVIATIONS FOR TANDEM AND PRIORITY QUEUES WITH GAUSSIAN INPUTS

By Michel Mandjes[1] and Miranda van Uitert

*CWI and University of Twente, and CWI and Vrije Universiteit*

This paper considers Gaussian flows multiplexed in a queueing network. A single node being a useful but often incomplete setting, we examine more advanced models. We focus on a (two-node) tandem queue, fed by a large number of Gaussian inputs. With service rates and buffer sizes at both nodes scaled appropriately, Schilder's sample-path large-deviations theorem can be applied to calculate the asymptotics of the overflow probability of the second queue. More specifically, we derive a lower bound on the exponential decay rate of this overflow probability and present an explicit condition for the lower bound to match the exact decay rate. Examples show that this condition holds for a broad range of frequently used Gaussian inputs. The last part of the paper concentrates on a model for a single node, equipped with a priority scheduling policy. We show that the analysis of the tandem queue directly carries over to this priority queueing system.

**1. Introduction.** Traffic engineering in communication networks greatly benefits from models that are capable of accurately describing and predicting the performance of the system. This modeling is a challenging task, as a broad variety of traffic types are multiplexed in the network, with each of them having its specific (stochastic) characteristics. A commonly used modeling step is to represent the network nodes as *queues*, and to use queueing theory to analyze the performance (in terms of loss, delay, throughput, etc.) of the nodes. For the single queue operating under the first-in-first-out (FIFO) discipline, even for advanced traffic models detailed analyses are available. Evidently this single-node FIFO model gives valuable insights, but is an oversimplification of reality. We mention two serious limitations.

[1]Part of this work was done while the first author was at Bell Laboratories/Lucent Technologies.

First, traffic streams usually traverse *concatenations* of hops (rather than just a single node). Second, it is envisaged that the service at these hops distinguishes between several traffic classes (by using *priority* mechanisms, or the more advanced *generalized processor sharing* discipline); compare the Differentiated Services (diffserv) approach proposed by the Internet Engineering Task Force [17]. This motivates the recent interest in performance evaluation for these more complex queueing models.

As indicated above, each type of traffic has its own stochastic properties, often summarized by the *correlation structure.* Traditional traffic models allow only *short-range dependent* traffic processes, such as Markov-modulated Poisson processes or exponential on-off sources, in which correlations decay relatively quickly. Traffic measurements in the 1990s, however, showed that in various situations *long-range dependent* traffic models are more appropriate. This explains the popularity of Gaussian models, as they cover both short-range (cf. Ornstein–Uhlenbeck) and long-range dependent models (e.g., fractional Brownian motion, see [19]). Another complicating issue is the fact that network traffic is usually influenced by feedback loops (think of TCP), which control how the user's traffic supply is transmitted into the network. Kilpi and Norros [18], however, argue that (nonfeedback) Gaussian traffic models are justified as long as the aggregation is sufficiently large (both in time and number of flows), due to central limit type of arguments.

This paper concentrates on the evaluation of tail asymptotics in queueing systems that are more advanced than a single FIFO node. More specifically, we examine in detail *tandem* queues (particularly the second queue) and *priority* queues (particularly the low-priority queue); it turns out that the analysis of the tandem queue essentially carries over to the priority system. Our paper is meant as a first step towards the analysis of networks with general topology, with nodes operating under advanced scheduling disciplines such as Generalized Processor Sharing (GPS).

In the tandem model we assume that $n$ i.i.d. Gaussian sources feed into the queueing system, where the (deterministic) service rates of the queues as well as the buffer thresholds are scaled by $n$, too. We now let $n$ go to infinity; the resulting framework is often referred to as the *many-sources* scaling, as was introduced in [32].

A vast body of results exists for single FIFO queues under the many-sources scaling. Most notably, under very mild conditions on the source behavior, it is possible to calculate the *exponential* decay of the probability $p_n(b,c)$ that the queue (fed by $n$ sources, and emptied at a deterministic rate $nc$) exceeds level $nb$. Early references in this large-deviations framework are the logarithmic asymptotics found in, for example, [7] and [8]. We remark that exact asymptotics for Gaussian inputs were recently found by Dębicki and Mandjes [9]. For Gaussian sources the logarithmic asymptotics of [7]

read

$$\lim_{n\to\infty}\frac{1}{n}\log p_n(b,c) = -\inf_{t>0}\frac{(b+(c-\mu)t)^2}{2v(t)}, \tag{1}$$

where $\mu$ is the mean input rate per source, and $v(t)$ is the variance of the amount of traffic generated by a single source in a time interval of length $t$. The goal of the present paper is to find expressions similar to (1) for tandem and priority queues.

Our work fits in the framework of a series of articles by Mannersalo and Norros [1, 23, 24, 25]. These papers examine queues with Gaussian sources, such as the single-node FIFO queue, but also priority queues and queues with GPS scheduling. For the latter types of queues, they derive heuristics for the decay rate of the overflow probabilities. The present paper shows that, for priority queues, these heuristics are typically close, but that there is a gap with the exact outcome. For both the tandem and priority queue a lower bound on the decay rate of the overflow probability is derived. In addition, we present an explicit condition under which this lower bound matches the exact value of the decay rate. Notice that lower bounds of the decay rate are usually of practical interest, as typically the network has to be designed such that overflow is sufficiently rare.

Our analysis exploits the above-mentioned similarity between priority and tandem queues. The techniques applied stem from large-deviations theory, particularly sample-path large deviations, based on (the generalized version of) Schilder's theorem. We mention that for priority systems in discrete time, different bounds were found by Wischik [33]; we will comment on the relation with our results later.

The paper is organized as follows. Section 2 introduces the tandem model, and presents preliminaries on (sample-path) large deviations. Section 3 analyzes the decay rate of the overflow probability of the second queue in a tandem system. This analysis is illustrated in Section 4 by a number of (analytical and numerical) examples. Section 5 studies the priority system, addressing the decay rate of the overflow probability in the low-priority queue.

## 2. Model and preliminaries.

This section introduces the tandem model that is analyzed in Section 3. In addition, we present preliminaries on large-deviations theory and the many-sources scaling.

### 2.1. *Tandem model.*

Consider a two-queue tandem model, with (deterministic) service rate $nc_1$ for the first queue and $nc_2$ for the second queue. We assume that $c_1 > c_2$, in order to exclude the trivial case where the buffer of the second queue cannot build up.

We consider $n$ sources (whose characteristics are specified in Section 2.2) that feed into the first queue. Traffic of these sources that has been served at the first queue immediately flows into the second queue—we assume no additional sources to feed the second queue. We are interested in the steady-state probability of the buffer content of the second queue $Q_{2,n}$ exceeding a certain threshold $nb$, $b > 0$, when the number of sources gets large, or, more specifically, its logarithmic asymptotics:

$$J := -\lim_{n\to\infty} \frac{1}{n} \log \mathbb{P}(Q_{2,n} > nb). \tag{2}$$

Note that we assume the buffer sizes of both queues to be infinite. We remark that it is not a priori clear that the limit in (2) exists; its existence is a result of our study (Theorem 3.1).

2.2. *Gaussian sources.* Let $A_i(\cdot)$ denote i.i.d. centered Gaussian processes with continuous sample paths and stationary increments, and $A_i(0) \equiv 0$, for $i = 1, \dots, n$. Then, for $s < t$, we interpret $A_i(s,t) := A_i(t) - A_i(s)$ as the amount of traffic generated by the $i$th source in $(s,t]$. Denote by $A(s,t)$ the generic random variable corresponding to a single source. The Gaussian sources are characterized by their *variance function* $v(\cdot)$ (which is necessarily continuous); for $s < t$, $\operatorname{Var} A(s,t) = v(t-s)$.

Although in this setup the Gaussian processes are centered, our analysis is capable of handling the situation in which the sources have a positive mean traffic rate $\mu$ (smaller than both $c_1$ and $c_2$, to guarantee stability). This is due to the fact that the results for centered sources can be translated immediately into results for noncentered sources; see Remark 2.6.

In the sequel we will frequently use the bivariate random variable $(A(t), A(s))$. It obviously obeys a two-dimensional Normal distribution with zero mean and covariance matrix $\Sigma(s,t)$. With $\Gamma(s,t) := \operatorname{Cov}(A(t), A(s))$, this covariance matrix is given by

$$\Sigma(s,t) := \begin{pmatrix} v(t) & \Gamma(s,t) \\ \Gamma(s,t) & v(s) \end{pmatrix} \quad \text{and} \quad \Gamma(s,t) = \frac{v(t) - v(|t-s|) + v(s)}{2}.$$

Gaussian sources have the conceptual problem that the possibility of *negative traffic* is not ruled out, as opposed to "classical" input processes, such as (compound) Poisson processes or on-off sources. However, in queueing theory a key role is played by functionals of the arrival process, which are well defined, regardless of whether the input stream corresponds to nonnegative traffic or not. Consider, for instance, the stationary distribution of a queue fed by a single source, emptied at rate $c$, given by the well-known formula $\sup_{t>0}(A(-t,0) - ct)$. Clearly, the distribution of such functionals can still be evaluated for Gaussian input; see, for example, Norros' pioneering work for fBm [26], or [15]. We remark that such an approach leads to nonnegative queue lengths in tandem systems with Gaussian inputs—this will follow directly from representation (10). For priority systems it is explained in detail in [25], Section 2.3, how negative queue lengths can be avoided (a discrete-time version of the priority discipline is introduced, in which negative traffic can annihilate queued traffic).

2.3. *Sample-path large deviations.* The analysis in the next sections relies on a sample-path large-deviations principle (LDP) for centered Gaussian processes. This section is devoted to a brief description of the main theorem in this field, (the generalized version of) *Schilder's theorem* [5]. However, we start by recalling (the multivariate version of) the well-known *Cramér's theorem*; see [10], Theorem 2.2.30. We let $\langle\cdot,\cdot\rangle$ denote the usual inner product: $\langle a,b\rangle := a^{\mathrm{T}}b = \sum_{i=1}^d a_i b_i$.

THEOREM 2.1 (Multivariate Cramér). *Let $X_i\in\mathbb{R}^d$ be i.i.d. $d$-dimensional random vectors, distributed as a random vector $X$ with finite moment-generating function $\mathbb{E}e^{\langle\theta,X\rangle}$ (for all $\theta\in\mathbb{R}^d$). Then the following LDP applies:*

(a) *For any closed set $F\subset\mathbb{R}^d$,*

$$\limsup_{n\to\infty}\frac{1}{n}\log\mathbb{P}\left(\frac{1}{n}\sum_{i=1}^n X_i\in F\right)\le-\inf_{x\in F}\Lambda(x).$$

(b) *For any open set $G\subset\mathbb{R}^d$,*

$$\liminf_{n\to\infty}\frac{1}{n}\log\mathbb{P}\left(\frac{1}{n}\sum_{i=1}^n X_i\in G\right)\ge-\inf_{x\in G}\Lambda(x),$$

*where the large deviations rate function $\Lambda(\cdot)$ is given by*

$$\Lambda(x):=\sup_{\theta\in\mathbb{R}^d}(\langle\theta,x\rangle-\log\mathbb{E}e^{\langle\theta,X\rangle}). \tag{3}$$

REMARK 2.2. Consider the specific case that $X$ has a multivariate Normal distribution with mean vector $\mu$ and $(d\times d)$ nonsingular covariance matrix $\Sigma$. Using $\log\mathbb{E}e^{\langle\theta,X\rangle}=\langle\theta,\mu\rangle+\frac{1}{2}\theta^{\mathrm{T}}\Sigma\theta$, it is not hard to derive that, with $(x-\mu)^{\mathrm{T}}\equiv(x_1-\mu_1,\dots,x_d-\mu_d)$,

$$\theta^\star=\Sigma^{-1}(x-\mu)\quad\text{and}\quad\Lambda(x)=\tfrac{1}{2}(x-\mu)^{\mathrm{T}}\Sigma^{-1}(x-\mu), \tag{4}$$

where $\theta^\star$ optimizes (3); it is well known that $\Lambda(\cdot)$ is convex.

We now sketch the framework of Schilder's sample-path LDP, as established in [5], see also [11]. We adopt the notation and setup of [1, 23].

Consider the $n$ i.i.d. centered Gaussian processes $A_i(\cdot)$, as introduced in Section 2.2. Define the path space $\Omega$ as

$$\Omega := \left\{\omega : \mathbb{R} \to \mathbb{R},\ \text{continuous},\ \omega(0) = 0, \lim_{t\to\infty} \frac{\omega(t)}{1+|t|} = \lim_{t\to-\infty} \frac{\omega(t)}{1+|t|} = 0\right\},$$

which is a separable Banach space by imposing the norm

$$\|\omega\|_\Omega := \sup_{t\in\mathbb{R}} \frac{|\omega(t)|}{1+|t|}.$$

In [1] it is pointed out that $A_i(\cdot)$ can be realized on $\Omega$ under the assumption that

$$\lim_{t\to\infty} \frac{v(t)}{t^\alpha} = 0 \qquad \text{for some } \alpha < 2. \tag{5}$$

We assume assumption (5) to be in force throughout this paper.

Next we introduce and define the *reproducing kernel Hilbert space* $R \subseteq \Omega$—see [3] for a more detailed account—with the property that its elements are roughly as smooth as the covariance function $\Gamma(s,\cdot)$. We start from a "smaller" space $R^\star$, defined by

$$R^\star := \left\{\omega : \mathbb{R} \to \mathbb{R}, \omega(\cdot) = \sum_{i=1}^n a_i \Gamma(s_i,\cdot), a_i, s_i \in \mathbb{R}, n \in \mathbb{N}\right\}.$$

The inner product on this space $R^\star$ is, for $\omega_a, \omega_b \in R^\star$, defined as

$$\langle \omega_a, \omega_b\rangle_R := \left\langle \sum_{i=1}^n a_i\Gamma(s_i,\cdot), \sum_{j=1}^n b_j\Gamma(s_j,\cdot)\right\rangle_R = \sum_{i=1}^n\sum_{j=1}^n a_i b_j \Gamma(s_i,s_j); \tag{6}$$

notice that this implies $\langle \Gamma(s,\cdot),\Gamma(\cdot,t)\rangle_R = \Gamma(s,t)$. This inner product has the following useful property, which we refer to as the *reproducing kernel* property:

$$\omega(t) = \sum_{i=1}^n a_i\Gamma(s_i,t) = \left\langle \sum_{i=1}^n a_i\Gamma(s_i,\cdot), \Gamma(t,\cdot)\right\rangle_R = \langle \omega(\cdot), \Gamma(t,\cdot)\rangle_R. \tag{7}$$

From this we introduce the norm $\|\omega\|_R := \sqrt{\langle\omega,\omega\rangle_R}$. The closure of $R^\star$ under this norm is defined as the space $R$. Now we can define the rate function of the sample-path LDP:

$$I(\omega) := \begin{cases} \frac{1}{2}\|\omega\|_R^2, & \text{if } \omega \in R, \\ \infty, & \text{otherwise.}\end{cases} \tag{8}$$

For the Gaussian sources introduced in Section 2.2, the following sample-path LDP holds.

THEOREM 2.3 (Generalized Schilder). *The following sample-path LDP applies:*

(a) *For any closed set* $F \subset \Omega$,

$$\limsup_{n\to\infty} \frac{1}{n} \log \mathbb{P}\left(\frac{1}{n}\sum_{i=1}^{n} A_i(\cdot) \in F\right) \le -\inf_{\omega\in F} I(\omega).$$

(b) *For any open set* $G \subset \Omega$,

$$\liminf_{n\to\infty} \frac{1}{n} \log \mathbb{P}\left(\frac{1}{n}\sum_{i=1}^{n} A_i(\cdot) \in G\right) \ge -\inf_{\omega\in G} I(\omega).$$

A difficulty of Schilder's theorem is its "implicitness," as only in special cases the rate function $I(\cdot)$ can be explicitly minimized over the set of interest. The authors of [1] succeed in exploiting the reproducing kernel property to give a sample-path analysis of overflow in a single FIFO queue (with deterministic service rate $nc$) fed by Gaussian inputs. With $Q_n$ denoting the stationary buffer content, they derive

$$\lim_{n\to\infty} \frac{1}{n} \log \mathbb{P}(Q_n > nb) = -\inf_{t\ge 0} \frac{(b+ct)^2}{2v(t)}.$$

It is elementary to show the existence of a minimizing $t$. On the one hand, it holds that $\lim_{t\downarrow 0}(b+ct)^2/2v(t) = \infty$ [due to $v(0) = 0$]. On the other hand, (5) implies that there is a $\beta < 2$ such that $v(t) < t^\beta$ eventually, and hence also $\lim_{t\to\infty}(b+ct)^2/2v(t) = \infty$. Notice that even in the FIFO setting, there is not necessarily uniqueness of the optimizing $t$; see, for instance, [21], [1], Section 3.7, or [20], Example 5.2.

If $t^\star$ denotes a minimizing $t$, the corresponding path is

$$f^\star(r) = -\frac{\Gamma(-r,t^\star)}{v(t^\star)}(b+ct^\star) = -\frac{v(t^\star) - v(|t^\star + r|) + v(-r)}{2v(t^\star)}(b+ct^\star). \tag{9}$$

This path corresponds to a buffer that starts to fill at time $-t^\star$, and reaches overflow at time 0; it is not hard to check that $f^\star(-t^\star) = -b - ct^\star$, and $f^\star(0) = 0$, as desired. Notice that the path is in $R$ (in fact even in $R^\star$). If there is a unique optimizing path in the target set (i.e., the set of all paths leading to overflow), it is usually referred to as the *most likely path to overflow*. It has the interpretation that, given that the rare event of overflow happens, with high probability it happens according to this trajectory. Also, $t^\star$ has then the interpretation of the most likely duration of the busy period preceding overflow. (Notice that, in this FIFO setting, there is not necessarily uniqueness; see, e.g., Section 3.7 in [1] or Example 5.2 in [20], or, in a non-Gaussian setting, [21].)

2.4. *Many-sources scaling.* In this section we show that the probability of our interest can be written in terms of the "empirical mean process" $n^{-1}\sum_{i=1}^n A_i(\cdot)$. The following lemma exploits the fact that we know both a representation of the first queue $Q_{1,n}$ (in steady state) and a representation of the *total* queue $Q_{1,n}+Q_{2,n}$ (in steady state). Let $t_0 := b/(c_1-c_2)$.

LEMMA 2.4. *$\mathbb{P}(Q_{2,n}>nb)$ equals*

$$\mathbb{P}\left(\exists\, t>t_0 : \forall\, s\in(0,t) : \frac{1}{n}\sum_{i=1}^n A_i(-t,-s) > b + c_2 t - c_1 s\right).$$

PROOF. Notice that a "reduction principle" applies: the total queue length is unchanged when the tandem network is replaced by its slowest queue; see [4, 14]. More formally: $Q_{1,n}+Q_{2,n}=\sup_{t>0}(\sum_{i=1}^n A_i(-t,0) - nc_2 t)$. Consequently we can rewrite

$$\begin{aligned} Q_{2,n} &= (Q_{1,n}+Q_{2,n}) - Q_{1,n} \\ &= \sup_{t>0}\left(\sum_{i=1}^n A_i(-t,0) - nc_2 t\right) - \sup_{s>0}\left(\sum_{i=1}^n A_i(-s,0) - nc_1 s\right). \end{aligned} \tag{10}$$

It was shown (see Lemma 5.1 of [31]) that the negative of the optimizing $t$ ($s$) corresponds to the start of the last busy period of the total queue (the first queue) in which time 0 is contained. Notice that a positive first queue induces a positive total queue, which immediately implies that we can restrict ourselves to $s\in(0,t)$. Hence $\mathbb{P}(Q_{2,n}>nb)$ equals

$$\mathbb{P}\left(\exists\, t>0 : \forall\, s\in(0,t) : \frac{1}{n}\sum_{i=1}^n A_i(-t,-s) > b + c_2 t - c_1 s\right).$$

Because for $s\uparrow t$ the requirement

$$\frac{1}{n}\sum_{i=1}^n A_i(-t,-s) > b + c_2 t - c_1 s$$

reads $0 > b + (c_2-c_1)t$, we can restrict ourselves to $t>t_0$. We can interpret $t_0$ as the minimum time it takes to cause overflow in the second queue (notice that the maximum net input rate of the second queue in a tandem system is $c_1-c_2$). $\square$

The crucial implication of the above lemma is that for analyzing $\mathbb{P}(Q_{2,n}\ge nb)$, we only have to focus on the behavior of the *empirical mean process*. More concretely,

$$\mathbb{P}(Q_{2,n}>nb) = \mathbb{P}\left(\frac{1}{n}\sum_{i=1}^n A_i(\cdot)\in\mathcal{S}\right), \tag{11}$$

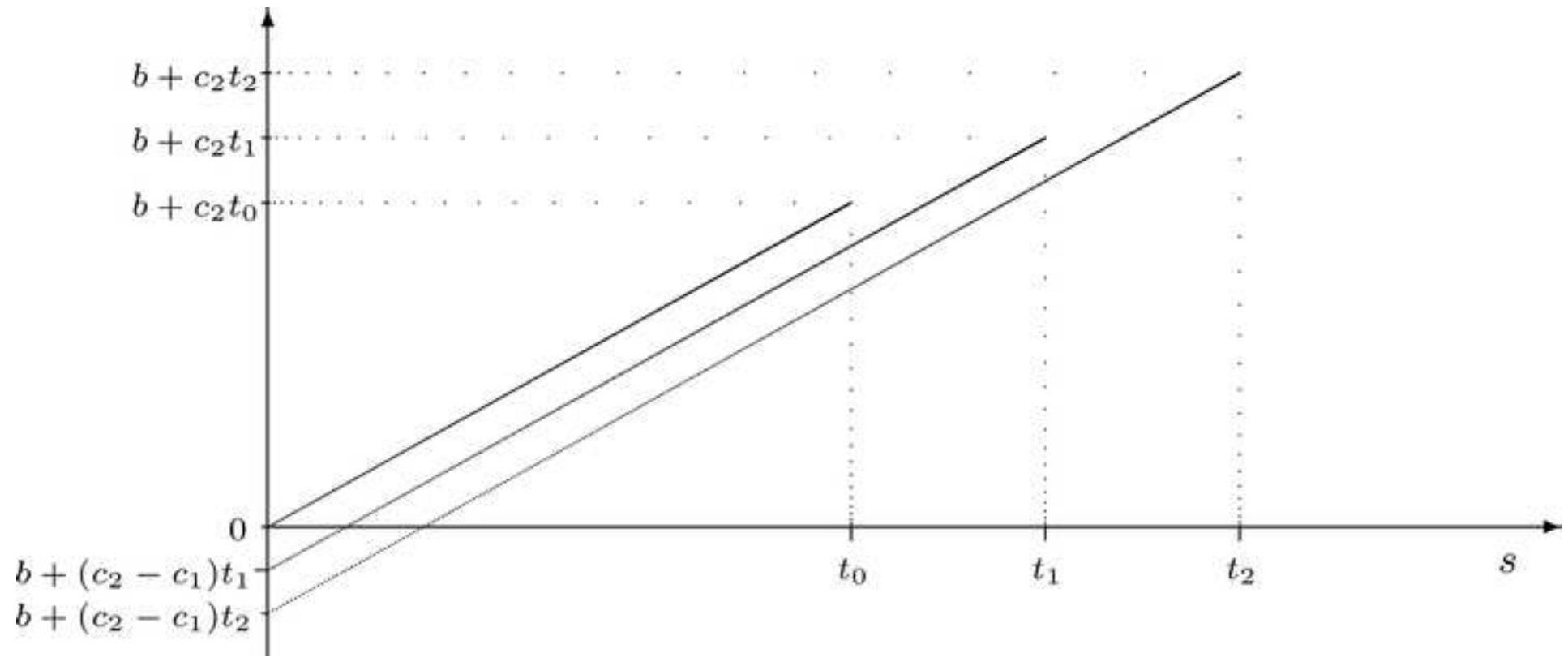

FIG. 1. *Graphical representation of the overflow set. For different values of $t$, the curve $b+c_2t-c_1(t-s)$ has been drawn. Overflow occurs if there is a $t>t_0$ such that the empirical mean process lies, for $s\in(0,t)$, above the corresponding curve.*

where the set of "overflow paths" $\mathcal{S}$ is given by

$$\mathcal{S}:=\{f\in\Omega:\exists\, t>t_0,\forall\, s\in(0,t):f(-s)-f(-t)>b+c_2t-c_1s\}.$$

REMARK 2.5. A straightforward time-shift shows that the probability that the empirical mean process is in $\mathcal{S}$ coincides with the probability that it is in $\mathcal{T}$, with

$$\mathcal{T}:=\{f\in\Omega:\exists\, t>t_0,\forall\, s\in(0,t):f(s)>b+c_2t-c_1(t-s)\}. \tag{12}$$

However, the set $\mathcal{T}$ is somewhat easier to interpret, see Figure 1. For different values of $t$ [i.e., $t_2>t_1>t_0=b/(c_1-c_2)$], the line $b+c_2t-c_1(t-s)$ has been drawn. The empirical mean process $n^{-1}\sum_{i=1}^n A_i(\cdot)$ is in $\mathcal{T}$ if there is a $t>t_0$ such that for all $s\in(0,t)$ it stays above the line $b+c_2t-c_1(t-s)$. Notice that $\mathcal{T}$ resembles the set corresponding to the probability of long busy periods in a single queue, as studied in [27].

REMARK 2.6. As indicated above, our results are for centered sources, but they can be translated easily into results for noncentered sources. Then the traffic generated by Gaussian source $i$ in the interval $[s,t)$ is $A(s,t)+\mu(t-s)$, where $A(s,t)$ corresponds to a centered source; here $0<\mu<\min\{c_1,c_2\}$ and $s<t$. Let $q(\mu,c_1,c_2)$ be the probability that the second queue exceeds $nb$, given that input rate $\mu$ and service rates $c_1$ and $c_2$ are in force. From (10) it follows immediately that $q(\mu,c_1,c_2)=q(0,c_1-\mu,c_2-\mu)$, and hence we can restrict ourselves to centered sources.

**3. Analysis.** In this section we analyze the logarithmic asymptotics of $\mathbb{P}(Q_{2,n} > nb)$. In Section 3.1 we show that the decay rate in (2) exists, of which we derive a lower bound in Section 3.2. It turns out that this lower bound has an insightful interpretation, which is given in Section 3.3. Section 3.4 presents conditions under which the lower bound is *tight* (meaning that the decay rate and lower bound match). Finally, in Section 3.5 we prove and explain some properties of the most likely path that we found.

3.1. *Decay rate of the overflow probability.* In this section we establish the existence of the decay rate (2) of $\mathbb{P}(Q_{2,n} > nb)$. We already saw in (11) that $\mathbb{P}(Q_{2,n} > nb)$ can be rewritten as the probability that the empirical mean process is in $\mathcal{S}$ (which is an open subset of $\Omega$). The existence of the decay rate follows from Schilder's result (Theorem 2.3), by showing that $\mathcal{S}$ is an $I$-continuity set, that is, that the infima of $I(\cdot)$ over $\mathcal{S}$ and $\overline{\mathcal{S}}$ match.

THEOREM 3.1.

$$J = \inf_{f \in \overline{\mathcal{S}}} I(f) = \inf_{f \in \mathcal{S}} I(f).$$

The proof of Theorem 3.1 can be found in the Appendix.

3.2. *Lower bound on the decay rate.* The main result of this section is a tractable lower bound on $J$, which is given in Theorem 3.2. Observe that

$$\mathcal{S} = \bigcup_{t > t_0} \bigcap_{s \in (0,t)} \mathcal{S}^{s,t} \qquad \text{with } \mathcal{S}^{s,t} := \{f \in \Omega : f(-s) - f(-t) > b + c_2 t - c_1 s\}.$$

Hence we are interested in the decay rate of the union of intersections. The decay rate of a union of events is simply the minimum of the decay rates of the individual events. The decay rate of an intersection is not standard. In the next theorem we find a straightforward lower bound on this decay rate. Define

$$\mathcal{U}^{s,t} := \{f \in \Omega : -f(-t) \geq b + c_2 t; f(-s) - f(-t) \geq b + c_2 t - c_1 s\}.$$

THEOREM 3.2. *The following lower bound applies:*

$$J \geq \inf_{t > t_0} \sup_{s \in (0,t)} \inf_{f \in \mathcal{U}^{s,t}} I(f). \tag{13}$$

PROOF. Clearly,

$$J = \inf_{t > t_0} \inf_{f \in \bigcap_{s \in (0,t)} \mathcal{S}^{s,t}} I(f).$$

Now fix $t$ and consider the inner infimum. If $f(-s) - f(-t) > b + c_2 t - c_1 s$ for all $s \in (0,t)$, then also ($f$ is continuous) $f(-s) - f(-t) \geq b + c_2 t - c_1 s$ for all $s \in [0,t]$. Hence,

$$\bigcap_{s\in(0,t)} \mathcal{S}^{s,t} \subseteq \bigcap_{s\in[0,t]} \mathcal{U}^{s,t} \subseteq \mathcal{U}^{r,t}$$

for all $r \in (0,t)$, and consequently

$$\inf_{f\in\bigcap_{s\in(0,t)}\mathcal{S}^{s,t}} I(f) \geq \inf_{f\in\mathcal{U}^{r,t}} I(f).$$

Now take the supremum over $r$ in the right-hand side. $\square$

Theorem 3.2 contains an infimum over $f \in \mathcal{U}^{s,t}$. In the next lemma we show how this infimum can be computed. Recalling (4), the bivariate large-deviations rate function of

$$\left(\sum_{i=1}^{n} \frac{A_i(-t,0)}{n}; \sum_{i=1}^{n} \frac{A_i(-t,-s)}{n}\right)$$

is, for $y, z \in \mathbb{R}$ and $t > 0$, $s \in (0,t)$, given by $\Lambda(y,z) := \frac{1}{2}(y,z)\Sigma(t-s,t)^{-1}(y,z)^{\mathrm{T}}$. We also introduce the following quantity, which plays a key role in our analysis:

$$\begin{aligned} k(s,t) &:= \mathbb{E}(A(-s,0) | A(-t,0) = b + c_2 t) \\ &= \mathbb{E}(A(s)|A(t) = b + c_2 t) = \frac{\Gamma(s,t)}{v(t)}(b + c_2 t). \end{aligned} \tag{14}$$

Under the following assumption, the infimum over $\mathcal{U}^{s,t}$ can be simplified considerably. The same assumption will be useful when deriving tightness conditions in Section 3.4.

Assumption 3.3. *$\sqrt{v(\cdot)} \in \mathcal{C}^2([0,\infty))$ is strictly increasing and strictly concave.*

Lemma 3.4. *Under Assumption 3.3, for $t > t_0$ and $s \in (0,t)$,*

$$\inf_{f\in\mathcal{U}^{s,t}} I(f) = \Upsilon(s,t) := \begin{cases} \Lambda(b + c_2 t, b + c_2 t - c_1 s), & \text{if } k(s,t) > c_1 s, \\ (b + c_2 t)^2 / 2v(t), & \text{if } k(s,t) \leq c_1 s. \end{cases}$$

Proof. Observe that

$$\begin{aligned} &\mathbb{P}\left(\sum_{i=1}^{n} \frac{A_i(\cdot)}{n} \in \mathcal{U}^{s,t}\right) \\ &\quad = \mathbb{P}\left(\sum_{i=1}^{n} \frac{A_i(-t,0)}{n} \geq b + c_2 t; \sum_{i=1}^{n} \frac{A_i(-t,-s)}{n} \geq b + c_2 t - c_1 s\right). \end{aligned} \tag{15}$$

Hence we can use Theorem 2.1, yielding

$$\inf_{f\in\mathcal{U}^{s,t}} I(f) = \inf \Lambda(y,z),$$

where the last infimum is over $y \geq b + c_2 t$ and $z \geq b + c_2 t - c_1 s$. Using that $\Lambda(\cdot,\cdot)$ is convex, this problem can be solved in a standard manner. It is easily verified that the contour of $\Lambda$ that touches the line $y = b + c_2 t$ does so at $z$-value

$$z_0 := \frac{\Gamma(t-s,t)}{v(t)}(b + c_2 t);$$

also the contour that touches $z = b + c_2 t - c_1 s$ does so at $y$-value

$$y_0 := \frac{\Gamma(t-s,t)}{v(t-s)}(b + c_2 t - c_1 s).$$

We first show that it cannot be that $y_0 > b + c_2 t$, as follows. If $y_0 > b + c_2 t$, then the optimum would be attained at $(y_0, b + c_2 t - c_1 s)$. Straightforward computations, however, show that $y_0 > b + c_2 t$ would imply that [use $\Gamma(t,t-s) \leq \sqrt{v(t)v(t-s)}$ ]

$$(\sqrt{v(t)} - \sqrt{v(t-s)})(b + c_2 t) > \sqrt{v(t)}c_1 s. \tag{16}$$

This inequality is not fulfilled for $s = 0$ $(0 \not> 0)$ nor for $s = t$ $(b + c_2 t \not> c_1 t$ for $t > t_0)$. As the left-hand side of (16) is convex (in $s$) due to Assumption 3.3, whereas the right-hand is linear (in $s$), there is no $s \in (0,t)$ for which the inequality holds. Conclude that $y_0 > b + c_2 t$ can be ruled out. Two cases are left:

(A) Suppose $z_0 > b + c_2 t - c_1 s$, or, equivalently, $k(s,t) \leq c_1 s$. Then $(b + c_2 t, z_0)$ is optimal (see the left panel of Figure 2), with rate function $(b + c_2 t)^2/2v(t)$, independent of $s$.

(B) In the remaining case (where $y_0 \leq b + c_2 t$ and $z_0 \leq b + c_2 t - c_1 s$) the optimum is attained at $(b + c_2 t, b + c_2 t - c_1 s)$, that is, the "corner point"; see the right panel in Figure 2. This happens if $k(s,t) > c_1 s$, and gives the desired decay rate.

This proves the stated. As an aside we mention that if $k(s,t) = c_1 s$, then both regimes coincide: $\Lambda(b + c_2 t, b + c_2 t - c_1 s) = (b + c_2 t)^2/2v(t)$. □

COROLLARY 3.5. *Under Assumption 3.3, the following lower bound applies:*

$$J \geq \inf_{t>t_0} \sup_{s\in(0,t)} \Upsilon(s,t).$$

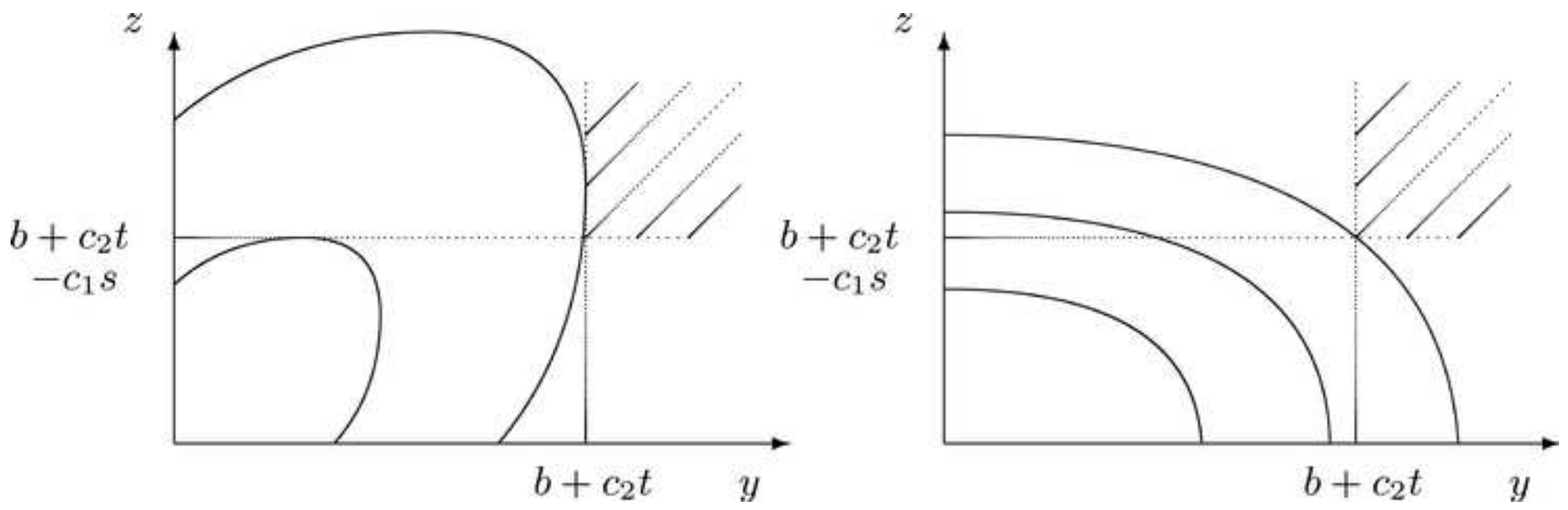

FIG. 2. *Contour lines of the (two-dimensional) rate function; the objective function is to be minimized over the shaded region.*

3.3. *Interpretation of the lower bound.* The results of the previous section have a helpful interpretation, leading to two regimes for values of $c_1$. For $c_1$ smaller than some critical link rate $c_1^F$, we show in Corollary 3.7 that the lower bound of Corollary 3.5 can be simplified considerably.

We start by drawing a parallel with the single-node FIFO result, as displayed in (1). There, $t$ has to be found such that

$$L_c(t) := \frac{(b + ct - \mathbb{E}A(t))^2}{2\,\mathrm{Var}\,A(t)}$$

is minimized. Let $t_c^F$ denote an optimizing argument $t$. $L_c(t)$ can be interpreted as the cost of generating $b + ct$ in an interval of length $t$, and $t_c^F$ as the time duration yielding the "lowest cost."

Now we turn to our tandem setting, and in particular to the result of Lemma 3.4. Computing the minimum of $\Lambda(y, z)$ over its admissible region, we saw that, under Assumption 3.3, in both cases the optimizing $y$ was equal to $y = b + c_2 t$. On the contrary, for the optimizing $z$ there were two possible regimes.

Now recall the representation (14) of $k(s, t)$ as a conditional mean, and (15). The result in Lemma 3.4 essentially states that in the regime $k(s, t) \le c_1 s$ the most likely realization of $\sum_{i=1}^n A_i(-t, 0) \ge nb + nc_2 t$ yields $\sum_{i=1}^n A_i(-t, -s) \ge nb + nc_2 t - nc_1 s$ (with high probability, $n$ large). In the other regime, $k(s, t) > c_1 s$, the most likely realization of $\sum_{i=1}^n A_i(-t, 0) \ge nb + nc_2 t$ does not automatically yield $\sum_{i=1}^n A(-t, -s) \le nb + nc_2 t - nc_1 s$ (with high probability, $n$ large); fulfilling the second constraint in (15) requires additional "cost."

The next decomposition result follows immediately from Lemma 3.4 and the above.

COROLLARY 3.6. *For $s \in (0, t)$, we have $\Upsilon(s, t) = L_{c_2}(t) + L(s|t)$, with*

$$L(s|t) := \frac{\max^2\{\mathbb{E}(A(s)|A(t) = b + c_2 t) - c_1 s, 0\}}{2\,\mathrm{Var}(A(s)|A(t) = b + c_2 t)}$$

$$= \frac{\max^2\{k(s,t) - c_1 s, 0\}}{2\,\mathrm{Var}(A(s)|A(t) = b + c_2 t)}. \tag{17}$$

Similarly to the interpretation of the single-node FIFO result, we can interpret $\Upsilon(s,t)$ as the cost of generating the required amount of traffic. Denoting by $s^\star$ and $t^\star$ optimizing arguments in Corollary 3.5, the intuition is as follows:

(A) "Cost component" $L_{c_2}(t)$ is needed to generate $b + c_2 t$ in the interval $(-t, 0]$. By taking the *infimum* over $t$ (to get $t^\star$) we find the *most likely* epoch to meet the constraint.

(B) "Cost component" $L(s|t)$ is required to make sure that no more than $c_1 s$ is generated in the interval $(-s, 0]$, *conditional* on the event $A(-t, 0) = b + c_2 t$. We can interpret $s^\star$ as the epoch at which *most* effort has to be done to fulfill this requirement. This is of course reflected by the fact that in Corollary 3.5 we have to take the *supremum* over all $s$ in $(0, t)$. Evidently, if $k(s,t) \le c_1 s$ for all $s \in (0,t)$, this cost component is 0.

For large values of $c_1$, $k(s,t)$ will be smaller than $c_1 s$ for all $s \in (0,t)$, since it does not depend on $c_1$. As argued above, in this case the second term in Corollary 3.6 vanishes. If this holds for the $t$ that maximizes the first term, that is, $t_{c_2}^F$, then

$$\inf_{t > t_0} \sup_{s \in (0,t)} \Upsilon(s,t) = L_{c_2}(t_{c_2}^F). \tag{18}$$

This clearly holds for all $c_1$ larger than

$$\begin{aligned} c_1^F &:= \inf\{c_1 | \forall\, s \in (0, t_{c_2}^F) : k(s, t_{c_2}^F) \le c_1 s\} \\ &= \inf\left\{c_1 | \forall\, s \in (0, t_{c_2}^F) : c_1 \ge \frac{k(s, t_{c_2}^F)}{s}\right\} = \sup_{s \in (0, t_{c_2}^F)} \frac{k(s, t_{c_2}^F)}{s}. \end{aligned}$$

It implies that, for these large values of $c_1$, the lower bound on $J$ of Corollary 3.5 coincides with the result of a single-node FIFO queue with service rate $c_2$. The intuition behind this is that essentially in this regime all traffic entering the first queue is served immediately, and goes directly into the second queue; traffic is not "reshaped" by the first queue. If $c_1 < c_1^F$, then the first queue *does* play a role in delaying and reshaping the traffic before entering the second queue, as we will see in the next section.

COROLLARY 3.7. *For all $c_1 \ge c_1^F$, (18) applies.*

3.4. *Tightness of the decay rate.* Corollary 3.5 is a lower bound on the decay rate $J$. Of course, such a bound is only useful if it is relatively close to the actual decay rate, or, even better, coincides with it. In the latter case we say that the lower bound is *tight*.

In Section 3.2, we have derived a lower bound on $J$ by replacing the decay rate of an intersection of events by the decay rate of the least likely of these. It is important to observe that if the optimum path in this least likely set happens to be in all the sets of the intersection, then the lower bound is tight.

More specifically, let $s^\star$ and $t^\star$ be optimizers in the lower bound of Corollary 3.5. Clearly we can prove tightness of the lower bound by showing that the most probable path in $\mathcal{U}^{s^\star,t^\star}$ is in $\mathcal{S}$ (or $\overline{\mathcal{S}}$; use Theorem 3.1). In our analysis we distinguish between (A) $c_1 \geq c_1^F$, and (B) $c_1 < c_1^F$.

Regime (A) ($c_1$ larger than the critical service rate). In this situation, we know from Corollary 3.7 that the lower bound in Corollary 3.5 reduces to the decay rate in a single FIFO queue. The next result follows easily.

Theorem 3.8. *Under Assumption 3.3, if $c_1 \geq c_1^F$, then*

$$J = \inf_{t>t_0} \sup_{s\in(0,t)} \Upsilon(s,t) = L_{c_2}(t_{c_2}^F),$$

*and a most probable path in $\mathcal{S}$ is*

$$f^\star(r) = -\mathbb{E}(A(r,0) | A(-t_{c_2}^F, 0) = b + c_2 t_{c_2}^F). \tag{19}$$

Proof. As shown in Section 3.3, in this regime $t^\star = t_{c_2}^F$, whereas the choice of $s^\star$ is irrelevant [as $c_1 \geq c_1^F$ implies $L(s|t^\star) = 0$ for all $s \in (0, t^\star)$]. Notice that it is now sufficient to show that $f^\star \in \mathcal{S}$, or $f^\star \in \overline{\mathcal{S}}$ (use Theorem 3.1). We claim that $f^\star(\cdot) \in \overline{\mathcal{S}}$, or more precisely, that there exists $t \geq t_0$ such that for all $s \in (0,t)$ it holds that $f^\star(-s) - f^\star(-t) \geq b + c_2 t - c_1 s$. This follows because, by definition of $c_1^F$, for all $s \in (0, t^\star)$,

$$\begin{aligned} f^\star(-s) - f^\star(-t^\star) &= \mathbb{E}(A(-t^\star, -s) | A(-t^\star, 0) = b + c_2 t^\star) \\ &= b + c_2 t^\star - k(s,t^\star) \geq b + c_2 t^\star - c_1 s. \end{aligned}$$

This completes the proof. □

We want to stress that the above theorem holds for all Gaussian processes, regardless of the specific shape of the variance function. Consequently, the result is also valid for long-range dependent processes, such as fractional Brownian motion.

REGIME (B) ($c_1$ smaller than the critical service rate). We follow the same approach as in Regime (A): first we derive (in Lemma 3.10) a most probable path in $\mathcal{U}^{s^\star,t^\star}$, and then we verify (in Theorem 3.11) whether this path is in $\mathcal{S}$. It turns out that we have to impose certain additional conditions to make the lower bound of Corollary 3.5 tight. We proceed by two technical lemmas; the proof of Lemma 3.9 is given in the Appendix.

LEMMA 3.9. *Under Assumption 3.3, if $c_1 < c_1^F$, then $k(s^\star, t^\star) \geq c_1 s^\star$.*

LEMMA 3.10. *If $k(s,t) \geq c_1 s$, then a most probable path in $\mathcal{U}^{s,t}$ is*

$$f(r) = -\mathbb{E}(A(r,0)|A(-t,0) = b + c_2 t, A(-s,0) = c_1 s), \tag{20}$$

*with norm $\Lambda(b + c_2 t, b + c_2 t - c_1 s)$.*

PROOF. Using standard properties of conditional multivariate Normal random variables, we see that $f(r)$ equals

$$\begin{aligned} &-\theta_1^\star(s,t)\Gamma(-r,t) - \theta_2^\star(s,t)\Gamma(-r,s) \\ &\qquad \text{with } \begin{pmatrix} \theta_1^\star(s,t) \\ \theta_2^\star(s,t) \end{pmatrix} := \Sigma(s,t)^{-1} \begin{pmatrix} b + c_2 t \\ c_1 s \end{pmatrix}. \end{aligned} \tag{21}$$

We finish the proof by applying Lemma 3.4, and observing that

$$\tfrac{1}{2}\|f\|_R^2 = \Upsilon(s,t) = \Lambda(b + c_2 t, b + c_2 t - c_1 s),$$

which is a matter of straightforward calculus. □

Before presenting our tightness result for the case $c_1 < c_1^F$, we introduce some new notation:

(i) For $r_1, r_2 < 0$,

$$\bar{\mathbb{E}} A(r_1, r_2) := \mathbb{E}(A(r_1, r_2)|A(-t^\star, 0) = b + c_2 t^\star),$$

with $\bar{\mathbb{V}\mathrm{ar}}(\cdot)$ and $\bar{\mathbb{C}\mathrm{ov}}(\cdot,\cdot)$ defined similarly. Also, $\bar{v}(r_1) := \bar{\mathbb{V}\mathrm{ar}}\, A(r_1, 0)$ and $\bar{\Gamma}(r_1, r_2) := \bar{\mathbb{C}\mathrm{ov}}(A(r_1,0), A(r_2,0))$.

(ii) For $r \in (-t^\star, 0)$ we define the functions

$$\bar{m}(r) := \frac{\bar{\mathbb{E}} A(r,0) + c_1 r}{\sqrt{\bar{v}(r)}}, \qquad m(r) := \frac{\bar{m}(r)}{\bar{m}(-s^\star)}, \qquad \rho(r) := \frac{\bar{\Gamma}(r, -s^\star)}{\sqrt{\bar{v}(r)\bar{v}(-s^\star)}}.$$

THEOREM 3.11. *Suppose*

$$m(-s) \leq \rho(-s) \qquad \text{for all } s \in (0, t^\star). \tag{22}$$

*Under Assumption 3.3, if $c_1 < c_1^F$, then*

$$J = \inf_{t > t_0} \sup_{s \in (0,t)} \Upsilon(s,t) = \Lambda(b + c_2 t^\star, b + c_2 t^\star - c_1 s^\star),$$

*and a most probable path is*

$$f^\star(r) = -\mathbb{E}(A(r,0)|A(-t^\star,0) = b + c_2 t^\star, A(-s^\star,0) = c_1 s^\star).$$

PROOF. As in Theorem 3.8, we have to show that $f^\star(\cdot)$ is in $\overline{\mathcal{S}}$. This is done as follows:

$$\begin{aligned} f^\star(-s) - f(-t^\star) &= \mathbb{E}(A(-t^\star,-s)|A(-t^\star,0) = b + c_2 t^\star, A(-s^\star,0) = c_1 s^\star) \\ &= b + c_2 t^\star - \bar{\mathbb{E}}(A(-s,0)|A(-s^\star,0) = c_1 s^\star) \\ &= b + c_2 t^\star - \bar{\mathbb{E}}A(-s,0) - \frac{\bar{\Gamma}(-s,-s^\star)}{\bar{v}(-s^\star)}(c_1 s^\star - \bar{\mathbb{E}}A(-s^\star,0)). \end{aligned}$$

Now it is easily seen that (22) implies that $f^\star(-s) - f(-t^\star) \geq b + c_2 t^\star - c_1 s$ for all $s \in (0, t^\star)$.

Due to Lemma 3.9, $k(s^\star, t^\star) \geq c_1 s^\star$. With Lemma 3.10, the expression for $J$ follows. □

Although the condition (22), required in Theorem 3.11, is stated in terms of the model parameters, as well as known statistics of the arrival process, it could be a tedious task to verify it in a specific situation. The next lemma presents a somewhat more transparent *necessary* condition for (22).

The intuition behind the lemma is the following. Observe that both $\rho(\cdot)$ and $m(\cdot)$ attain a maximum 1 at $r = -s^\star$. For $\rho(\cdot)$ this follows from the observation that $\rho(r)$ is a correlation coefficient; for $m(\cdot)$ from Corollary 3.6 and Lemma 3.9. Then a necessary condition for (22) is that in $s^\star$ the curve $m(\cdot)$ is "more concave" than $\rho(\cdot)$. The proof of the lemma is given in the Appendix.

LEMMA 3.12. *A necessary condition for* (22) *is*

$$m''(-s^\star) \leq \rho''(-s^\star), \tag{23}$$

*or equivalently,*

$$\theta_1^\star(s^\star,t^\star)(v''(t^\star - s^\star) - v''(s^\star)) + \theta_2^\star(s^\star,t^\star)(v''(0) - v''(s^\star)) \geq 0. \tag{24}$$

Condition (24) has an insightful interpretation, which will be given in the next section.

3.5. *Properties of the input rate path.* So far, we have analyzed paths $f$ of the *cumulative* amount of traffic injected into the system. In this section we turn our attention to the first derivative of $f$, which can be interpreted as the path of the *input rate* of the queueing system. As before, we have to consider two regimes: (A) $c_1 \geq c_1^F$, and (B) $c_1 < c_1^F$; let Assumption 3.3 be in force. Consider the paths $f^\star$ as identified in Theorems 3.8 and 3.11, and, more specifically, their derivative $g^\star(\cdot) := (f^\star)'(\cdot)$. In case (A), with $t^\star = t^F_{c_2}$, and $r \in (-t^\star, 0)$,

$$g^\star(r) = \frac{b + c_2 t^\star}{2v(t^\star)}(v'(r + t^\star) + v'(-r)),$$

whereas in case (B) it turns out that, with $r \in (-t^\star, -s^\star]$,

$$g^\star(r) = \frac{v'(r + t^\star) + v'(-r)}{2}\theta^\star_1(s^\star, t^\star) + \frac{-v'(-r - s^\star) + v'(-r)}{2}\theta^\star_2(s^\star, t^\star),$$

and with $r \in [-s^\star, 0)$,

$$g^\star(r) = \frac{v'(r + t^\star) + v'(-r)}{2}\theta^\star_1(s^\star, t^\star) + \frac{v'(r + s^\star) + v'(-r)}{2}\theta^\star_2(s^\star, t^\star).$$

If $v'(0) = 0$, we show below that the path $g^\star(\cdot)$ has some nice properties. Notice that the requirement $v'(0) = 0$ holds for many Gaussian processes. It is not valid for standard Brownian motion (Bm), since then $v(t) = t$, but the special structure of Bm allows an explicit analysis, see Section 4.1. Fractional Brownian motion (fBm), with $v(t) = t^{2H}$, has $v'(0) = 0$ only for $H \in (\frac{1}{2}, 1]$; see Section 4.2.

PROPOSITION 3.13. *If $c_1 \geq c^F_1$ and $v'(0) = 0$, then $g^\star(0) = g^\star(-t^\star) = c_2$.*

PROOF. Notice that, due to (1), $t^\star$ satisfies

$$2c_2 \frac{v(t^\star)}{v'(t^\star)} = b + c_2 t^\star.$$

The stated follows immediately from $v'(0) = 0$. [As an aside, we mention that $g^\star(\cdot)$ is symmetric in $-t^\star/2$.] □

Just as we exploited properties of $t^\star$ in the proof of Proposition 3.13, we need conditions for $s^\star$ and $t^\star$ in the regime $c_1 < c^F_1$. These are derived in the next lemma.

LEMMA 3.14. *If $c_1 < c^F_1$, then $s^\star$ and $t^\star$ satisfy the following equations:*

$$2c_2 = \theta^\star_1(s^\star, t^\star)v'(t^\star) + \theta^\star_2(s^\star, t^\star)(v'(t^\star) - v'(t^\star - s^\star)),$$
$$2c_1 = \theta^\star_2(s^\star, t^\star)v'(s^\star) + \theta_1(s^\star, t^\star)(v'(s^\star) + v'(t^\star - s^\star)).$$

PROOF. By Lemma 3.9, $k(s^\star, t^\star) \geq c_1 s^\star$. Observe that $\Upsilon(s, t) = \Lambda(b + c_2 t, b + c_2 t - c_1 s)$ can be rewritten as

$$\theta^{\mathrm{T}} x(s, t) - \tfrac{1}{2}\theta^{\mathrm{T}}\Sigma(s, t)\theta \qquad \text{where } x(s, t) := \begin{pmatrix} b + c_2 t \\ c_1 s \end{pmatrix}; \tag{25}$$

here we abbreviate $\theta \equiv (\theta_1^\star(s,t), \theta_2^\star(s,t))^{\mathrm{T}}$. We write $\partial_t$ and $\partial_s$ for the partial derivatives with respect to $t$ and $s$, respectively. The optimal $s^\star$ and $t^\star$ necessarily satisfy the first-order conditions, obtained by differentiating (25) to $t$ and $s$, and equating them to 0. Direct calculations yield

$$\begin{pmatrix} \theta_1 c_2 \\ \theta_2 c_1 \end{pmatrix} = \begin{pmatrix} \partial_t\theta_1 & \partial_t\theta_2 \\ \partial_s\theta_1 & \partial_s\theta_2 \end{pmatrix} (\Sigma(s,t)\theta - x(s,t)) + \begin{pmatrix} \frac{1}{2}\theta_1^2 v'(t) + \partial_t\Gamma(s,t)\theta_1\theta_2 \\ \frac{1}{2}\theta_2^2 v'(s) + \partial_s\Gamma(s,t)\theta_1\theta_2 \end{pmatrix}.$$

The second equality in (21) provides $x(s,t) = \Sigma(s,t)\theta$. Now the stated follows directly. □

PROPOSITION 3.15. *If* $c_1 < c_1^F$ *and* $v'(0) = 0$, *then* (i) $g^\star(-t^\star) = c_2$, *and* (ii) $g^\star(-s^\star) = c_1$. *Also, the necessary condition* (24) *is equivalent to* $(g^\star)'(-s^\star) \geq 0$.

PROOF. Claims (i) and (ii) follow directly from $v'(0) = 0$ and Lemma 3.14. The last statement follows directly after some calculations. □

Proposition 3.15 can be interpreted as follows. The second queue starts a busy period at time $-t^\star$. During this trajectory, the first queue starts to fill at time $-s^\star$ and is empty again at time 0, if the conditions of Theorem 3.11 apply. It is also easily seen that the necessary condition (24) has the appealing interpretation that $(g^\star)'(-s^\star) \geq 0$: the input rate path should be increasing at time $-s^\star$.

### 3.6. *Some remarks.*

REMARK 3.16. In our lower bound we replace the intersection over $s \in (0,t)$ by the *least likely event* of the intersection. Under condition (24) the occurrence of the least likely event *implies all the other events in the intersection*, with high probability [in the sense that $f^\star \in \mathcal{U}^{s^\star,t^\star}$ implies that $f^\star \in \mathcal{U}^{s,t^\star}$ for all $s \in (0,t^\star)$]. The examples in Section 4 show that (22) is met for many "standard" Gaussian models, but not always. If there is no tightness, a better lower bound can be obtained by approximating the intersection by more than just one event:

$$J \geq \inf_{t>t_0} \sup_{\overline{s}\in(0,t)^m} \inf_{f\in\mathcal{U}^{\overline{s},t}} I(f),$$

where $\overline{s} = (s_1, \ldots, s_m)$, and the "multiple-constraints set" $\mathcal{U}^{\overline{s},t}$ is defined by

$$\mathcal{U}^{\overline{s},t} := \{f \in \Omega : -f(-t) \geq b + c_2 t;$$
$$f(-s_i) - f(-t) \geq b + c_2 t - c_1 s_i, \text{ for } i = 1, \ldots, m\}.$$

Obviously, the lower bound becomes tighter when increasing $m$.

Remark 3.17. The approach we have followed in this section to analyze the two-node tandem network can be easily applied to an $m$-node tandem network, with strictly decreasing service rates, that is, $c_1 > \cdots > c_m$—nodes $i$ for which $c_i \leq c_{i+1}$ can be ignored, see [4, 14, 16]. Note that $\sum_{i=1}^{k} Q_{i,n}$ is equivalent to the FIFO queue in which the sources feed into a buffer that is emptied at rate $c_k$. This means that we have the characteristics of both $\sum_{i=1}^{m-1} Q_{i,n}$ and $\sum_{i=1}^{m} Q_{i,n}$, which enables the analysis of $Q_{m,n}$, just as in the two-node tandem case.

## 4. Examples.

One of the reasons for considering Gaussian input processes is that they cover a broad range of correlation structures. Choosing the variance function appropriately, we can make the input process exhibiting, for instance, long-range dependent behavior. In this section we do the computations for various variance functions. We also discuss in detail the condition in Theorem 3.11.

### 4.1. *Standard Brownian motion.*

The variance function for Brownian motion (Bm) is given by $v(t) = t$. Using (1), it is easily found that $t_{c_2}^F = b/c_2$. According to Corollary 3.7, $c_1^F$ is the largest value of $c_1$ such that for all $s \in (0, t_{c_2}^F)$,

$$\frac{s}{t_{c_2}^F}(b + c_2 t_{c_2}^F) - c_1 s \leq 0,$$

that is, $c_1^F = 2c_2$. Hence, using Theorem 3.8, we have for $c_1 \geq 2c_2$ that $J = 2bc_2$, with a constant input rate $g^\star(r) = 2c_2$ for $r \in (-t_{c_2}^F, 0)$ and $g^\star(r) = 0$ elsewhere.

Now we turn to the case where $c_1 < 2c_2$. The optimizing $s^\star$ and $t^\star$ are determined by solving the first-order equations for $s$ and $t$; see Theorem 3.11. We immediately obtain that $t^\star = b/(c_1 - c_2)$ and $s^\star = 0$. Obviously, for this regime the service rate of the first queue *does* play a role. The most probable input rate path reads $g^\star(r) = c_1$, for $r \in (-t^\star, 0)$ and $g^\star(r) = 0$ elsewhere. It is easily verified that the most probable path $f^\star(\cdot)$ is in $\overline{\mathcal{S}}$, making the decay rate as found in Theorem 3.11 tight. In other words,

$$J = \Lambda(b + c_2 t^\star, b + c_2 t^\star - c_1 s^\star) = \frac{bc_1^2}{2(c_1 - c_2)}.$$

Observe that, interestingly, Bm apparently changes its rate instantaneously, as reflected by the most likely input rate path. This is a consequence of the independence of the increments.

4.2. *Fractional Brownian motion.* The variance function for fractional Brownian motion (fBm) is given by $v(t) = t^{2H}$, where $H$ is the so-called Hurst parameter. For $H > \frac{1}{2}$ this corresponds to long-range dependent traffic. Now (1) gives

$$t_{c_2}^F = \frac{b}{c_2}\frac{H}{1-H}.$$

By Theorem 3.8,

$$J = \frac{1}{2}\left(\frac{b}{1-H}\right)^{2-2H}\left(\frac{c_2-\mu}{H}\right)^{2H}$$

for all $c_1 \geq c_1^F$. Unfortunately, for general $H$ there does not exist a closed-form expression for $c_1^F$. Now turn to the case $c_1 < c_1^F$. Lemma 3.12 states that (24) is a necessary condition for tightness to hold. Observe that $v''(t) = (2H-1)2Ht^{2H-2}$ and hence $v''(0) = \infty$. It is easily checked that $\theta_2^\star(s^\star, t^\star) \leq 0$, which implies that in this case (24) is not satisfied. Therefore the lower bound on $J$ is *not* tight.

4.3. M/G/$\infty$ *input.* A versatile traffic model is the so-called M/G/$\infty$ input process. In this model sessions arrive according to a Poisson process with rate $\lambda$, and stay in the system for some random duration $D$. During this period they generate traffic at a unit rate. By choosing specific session-length distributions $D$, both short-range and long-range dependent inputs can be modeled. For more results on queues with M/G/$\infty$ input traffic processes, see, for example, [12, 30]. Below we approximate the M/G/$\infty$ inputs by their "Gaussian counterpart," that is, Gaussian sources with the same mean and variance as the M/G/$\infty$ input; this procedure is extensively motivated in [1, 2].

Let the mean session-length be finite, say $\delta$, such that the mean input rate equals $\lambda\delta$. We denote by $F_D(\cdot)$ the distribution function of $D$ and by $F_{D^r}(\cdot)$ the distribution function of the *residual* session-length, that is, $F_{D^r}(x) = \delta^{-1}\int_0^x (1 - F_D(y))\,dy$. We denote the corresponding densities by $f_D(\cdot)$ and $f_{D^r}(\cdot)$.

Let $B(t)$ denote the amount of traffic generated by a single M/G/$\infty$ input in an interval of length $t$. We now show how to compute the variance $v(\cdot)$ of $B(t)$. We will do this by first deriving the moment-generating function of $B(t)$. In fact two types of sources contribute:

1. Sources that were already present at the start of the interval. The number of these sources has a Poisson distribution with mean $\lambda\delta$. Their residual duration has density $f_{D^r}(\cdot)$; with probability $(1 - F_{D^r}(t))$ they transmit traffic during the entire interval.

2. Sources that arrive during the interval. Their number has a Poisson($\lambda t$) distribution. Given that the number of these arrivals is $k \in \mathbb{N} \cup \{0\}$, their arrival epochs are i.i.d. random variables, uniformly over the interval (with density $t^{-1}$). Their duration has density $f_D(\cdot)$.

Straightforward computations now yield (cf. [22])

$$\log \mathbb{E}(e^{\theta B(t)}) = \lambda\delta(M_t(\theta) - 1) + \lambda t(N_t(\theta) - 1)$$

with

$$M_t(\theta) := \int_0^t e^{\theta x} f_{D^r}(x)\, dx + e^{\theta t}(1 - F_{D^r}(t))$$

and

$$N_t(\theta) := \int_0^t \int_u^t \frac{1}{t} e^{\theta(x-u)} f_D(x-u)\, dx\, du + \int_0^t \frac{1}{t} e^{\theta(t-u)}(1 - F_D(t-u))\, du.$$

Taking the second derivative of the log moment-generating function (with respect to $\theta$) and then substituting 0 for $\theta$, gives the variance $v(t)$ of $B(t)$:

$$\lambda\delta\biggl(\int_0^t x^2 f_{D^r}(x)\, dx + t^2(1 - F_{D^r}(t))\biggr)$$
$$+ \lambda\biggl(\int_0^t \int_u^t (x-u)^2 f_D(x-u)\, dx\, du + \int_0^t (t-u)^2(1 - F_D(t-u))\, du\biggr).$$

For fBm we could a priori rule out tightness of the lower bound due to $v''(0) = \infty$; see Lemma 3.12. For M/G/$\infty$ inputs we show in the following lemma that $v''(0)$ is finite, even for heavy-tailed $D$. It implies that condition (22) needs to be checked to verify tightness.

Lemma 4.1. *For $\delta < \infty$ and finite $f_D(\cdot)$, both* (i) $v'(0) = 0$ *and* (ii) $v''(0) < \infty$.

Proof. Using standard rules for differentiation of integrals,

$$v'(t) = \lambda\delta 2t(1 - F_{D^r}(t)) + \lambda \int_0^t 2(t-u)(1 - F_D(t-u))\, du$$

and hence $v'(0) = 0$. Similarly,

$$\begin{aligned} v''(t) &= 2\lambda\delta(1 - F_{D^r}(t) - t f_{D^r}(t)) \\ &\quad + 2\lambda \int_0^t (1 - F_D(t-u) - (t-u) f_D(t-u))\, du \\ &= 2\lambda \int_t^\infty (1 - F_D(s))\, ds. \end{aligned}$$

Hence, $v''(0) = 2\lambda\delta < \infty$. □

Now we consider some examples of session-length distributions. In all the examples we take $b = 0.5$, $\lambda = 0.125$, $\delta = 2$ and $c_2 = 1$.

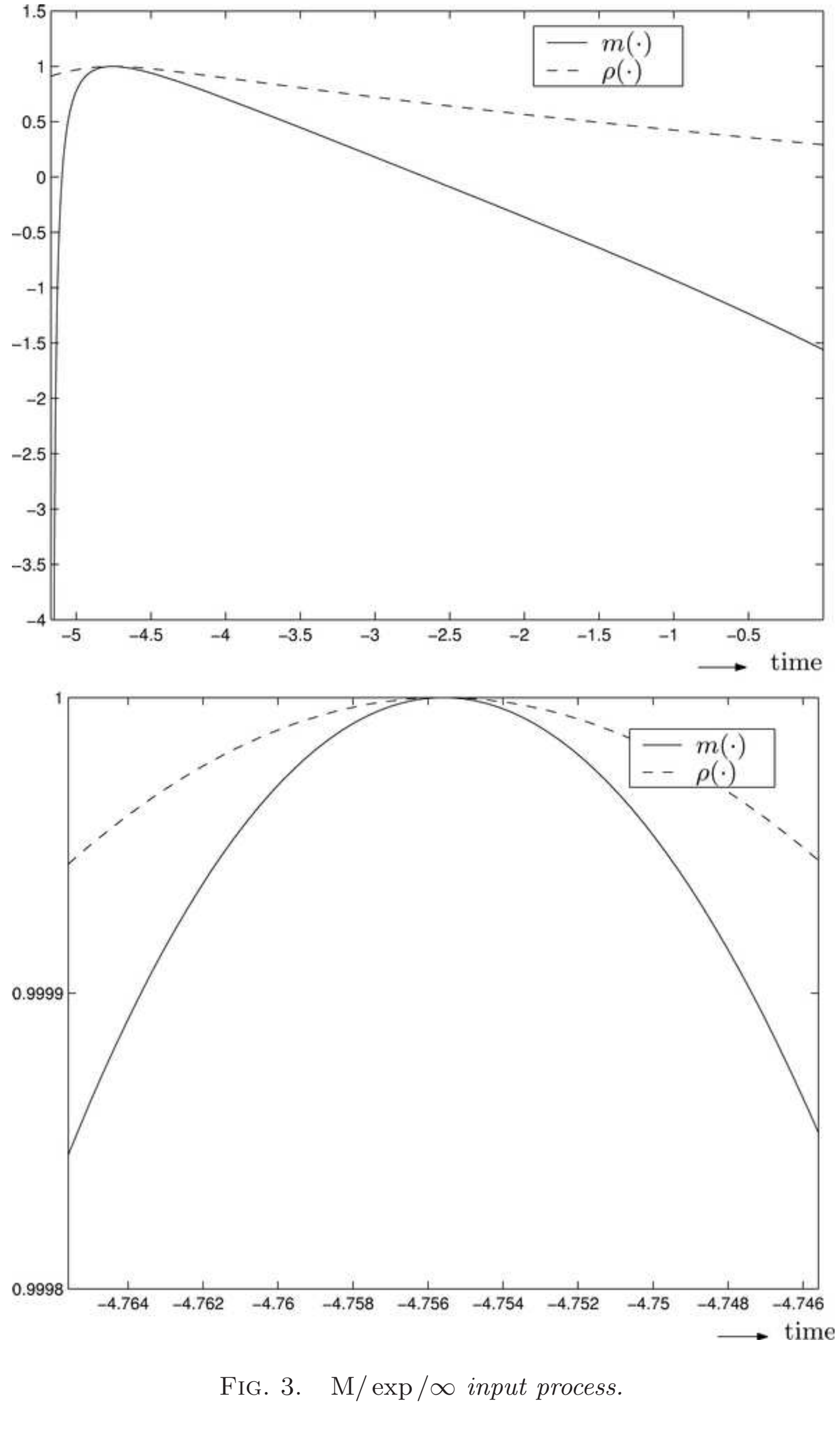

FIG. 3. M/exp/$\infty$ *input process.*

*Exponential.* Using the above formula for $v(\cdot)$, we get

$$v(t) = 2\lambda\delta^3\left(\frac{t}{\delta} - 1 + \exp\left(-\frac{t}{\delta}\right)\right).$$

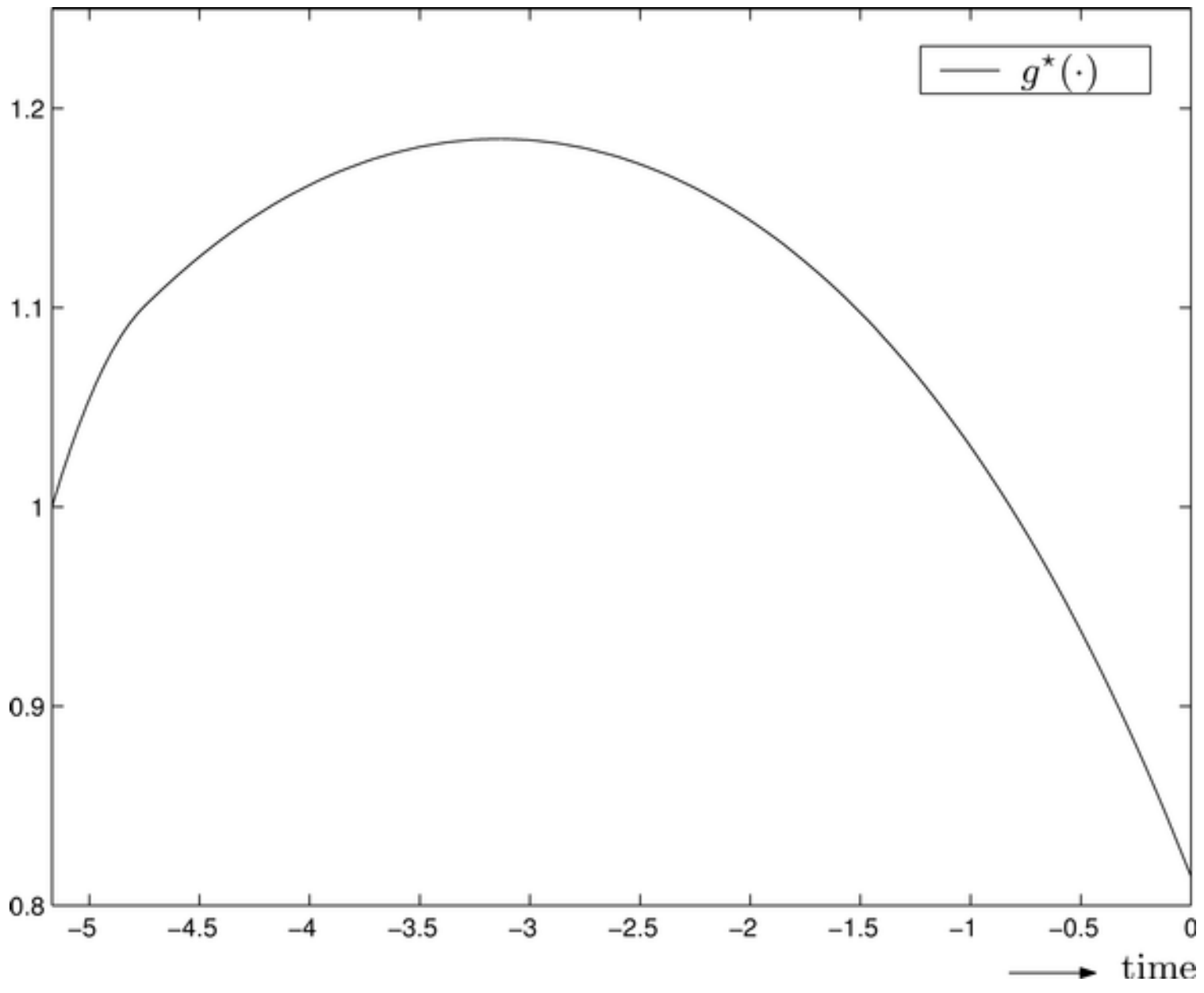

Fig. 4. *Input rate path for* M/ exp /∞ *input process.*

Notice that $v(\cdot)$ tends to a straight line for large $t$ (corresponding to short-range dependence). Numerical computations then give $c_1^F = 1.195$. Taking $c_1 = 1.1$ results in $s^\star = 4.756$, $t^\star = 5.169$ and $m(r)$, $\rho(r)$ as given in Figure 3. The upper panel of Figure 3 shows $m(r)$ and $\rho(r)$ for $r \in (-t^\star, 0)$, whereas the lower panel magnifies the graph around $-s^\star$. We see that indeed $m(\cdot) \leq \rho(\cdot)$ on the desired interval, so the decay rate is tight. A corresponding input rate path is given in Figure 4, which satisfies the properties as indicated in Proposition 3.15.

*Hyperexponential.* In case $D$ has a hyperexponential distribution, with probability $p_i \in (0, 1)$ it behaves as an exponential random variable with mean $\nu_i^{-1}$, with $i = 1, 2$ and $p_1 + p_2 = 1$. It is easily verified that

$$v(t) = 2\lambda \frac{p_1}{\nu_1^3}(\nu_1 t - 1 + e^{-\nu_1 t}) + 2\lambda \frac{p_2}{\nu_2^3}(\nu_2 t - 1 + e^{-\nu_2 t}),$$

with $\nu_2 = p_2/(\delta - p_1/\nu_1)$. As in the exponential case, $v(\cdot)$ is asymptotically linear. For $p_1 = 0.25$ and $\nu_1 = 5$, we find $c_1^F = 1.173$, and $s^\star = 4.700$, $t^\star = 5.210$, when using $c_1 = 1.1$. Also for this example $m(\cdot) \leq \rho(\cdot)$, and hence there is tightness; the graph looks similar to Figure 3. A corresponding input rate path is given in Figure 5.

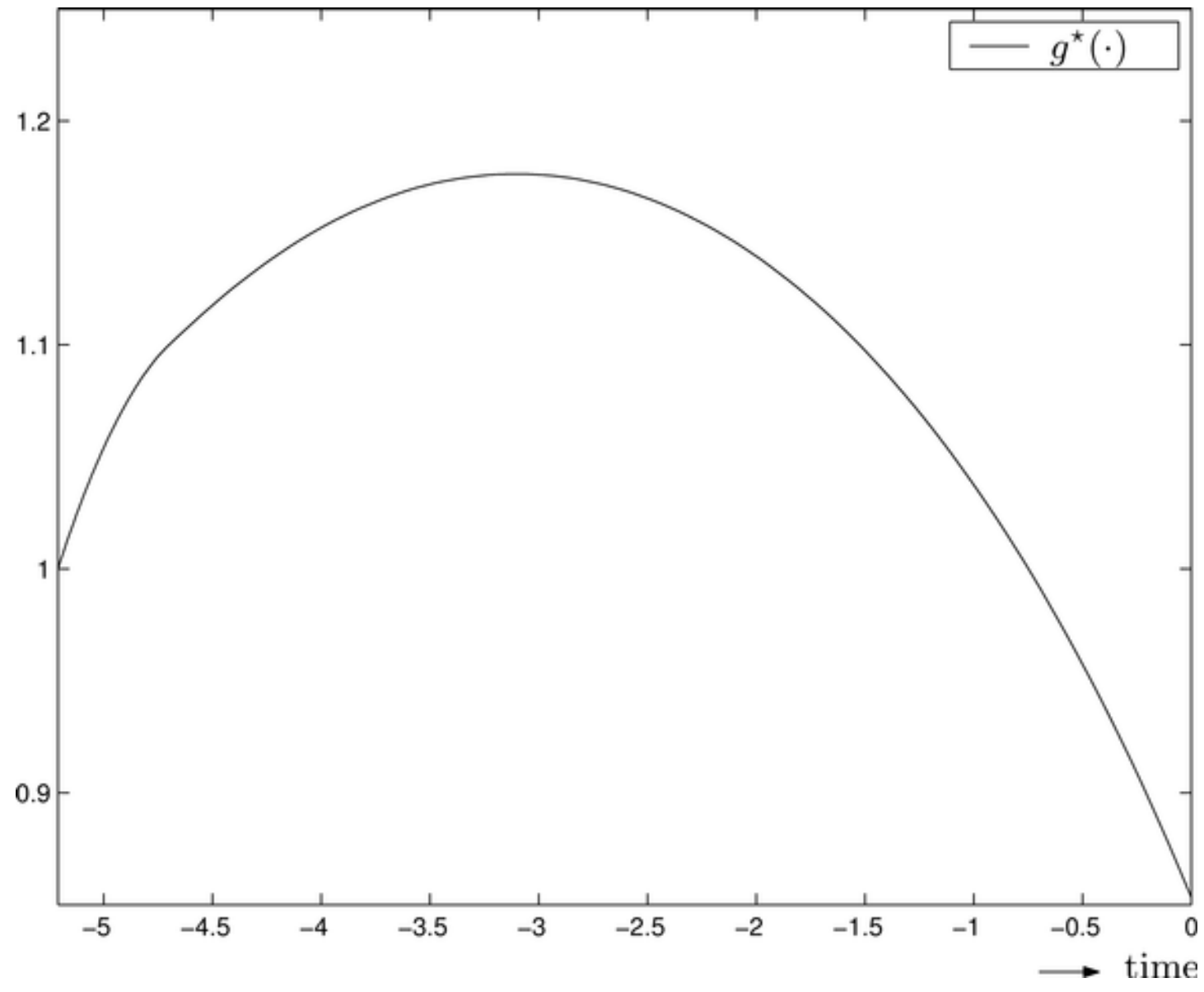

Fig. 5. *Input rate path for* M/H2/∞ *input process.*

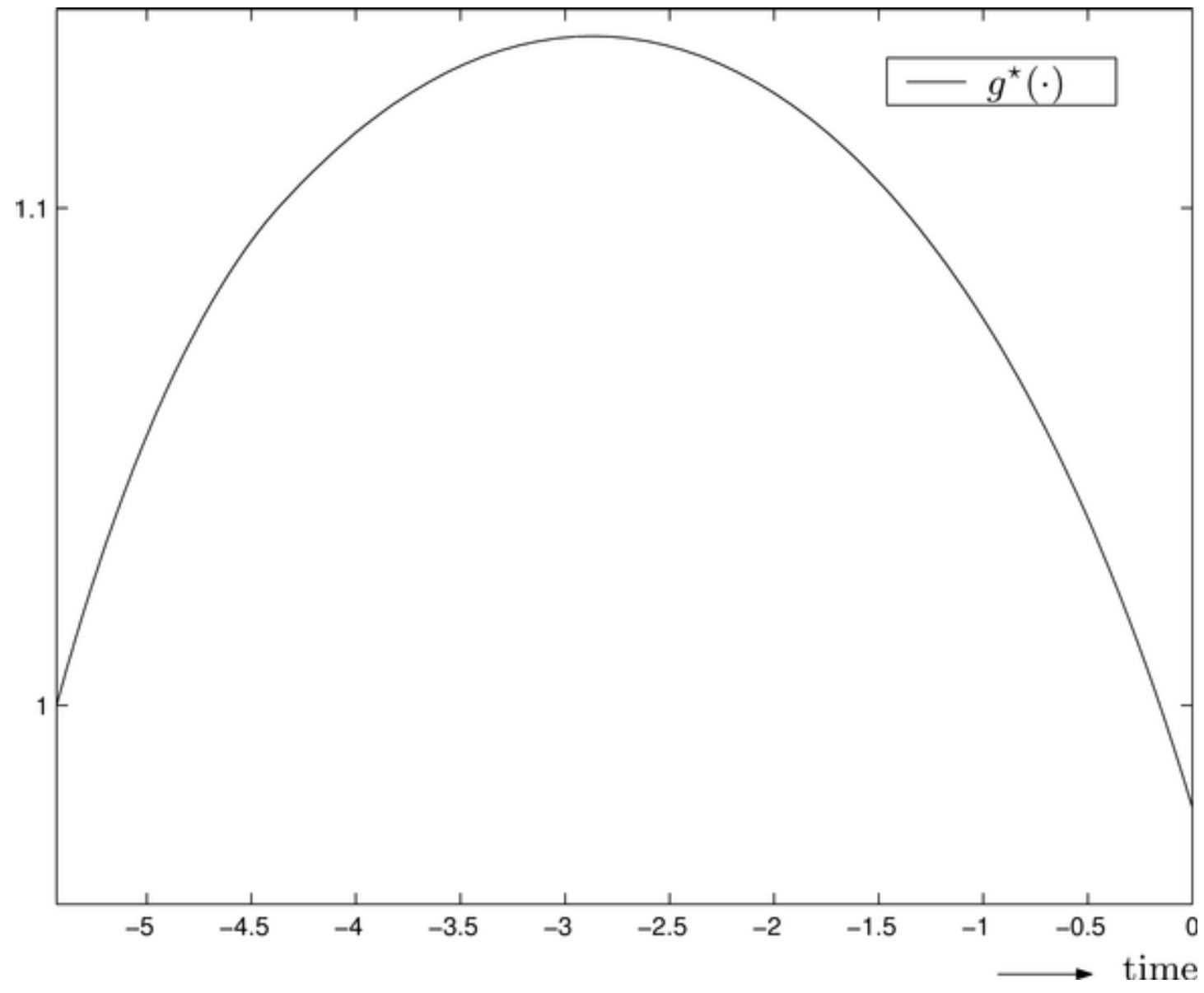

Fig. 6. *Input rate path for* M/Par/∞ *input process.*

*Pareto.* If $D$ has a Pareto distribution, then $\mathbb{P}(D > t) = (1/(1+t))^\alpha$. The variance function is given by

$$v(t) = \frac{2\lambda}{(3-\alpha)(2-\alpha)(1-\alpha)}(1-(t+1)^{3-\alpha}+(3-\alpha)t),$$

with $\alpha = (1+\delta)/\delta$, excluding $\delta = 1$ or $\frac{1}{2}$. Notice that we have $\alpha = 1\frac{1}{2}$, yielding $v(t) \sim t\sqrt{t}$, which corresponds to long-range dependent traffic. Numerical calculations show that $c_1^F = 1.115$, and for $c_1 = 1.1$ we obtain $s^\star = 4.373$, $t^\star = 5.432$. Again $m(\cdot)$ is majorized by $\rho(\cdot)$. An input rate path is given in Figure 6. We empirically found that there is not always tightness in the M/Par/$\infty$ case. If $b$ is larger, for instance $b = 1$, then (22) is not met.

**5. Priority queues.** In Section 3 we analyzed overflow in the second queue of a tandem system. This analysis was enabled by the fact that we had explicit knowledge of both the *first* queue and the *total* queue. In the present section we use the same type of arguments to solve the (two-queue) priority system.

5.1. *Analysis.* We consider a priority system with a link of capacity $nc$, fed by traffic of two classes, each with its own queue. Traffic of class 1 does not "see" class 2 at all, and consequently we know how the *high-priority* queue $Q_{h,n}$ behaves. Also, due to the work-conserving property of the system, the *total* queue length $Q_{h,n} + Q_{\ell,n}$ can be characterized. Now we are able, applying the same arguments as for the tandem queue, to analyze the decay rate of the probability of exceeding some buffer threshold in the low-priority queue. This similarity between tandem and priority systems has been observed before; see, for instance, [13].

We let the system be fed by $n$ i.i.d. high-priority (hp) sources, and an equal number of i.i.d. low-priority (lp) sources; both classes are independent. We assume that both hp and lp sources are Gaussian, and satisfy the requirements imposed in Section 2. Define the means by $\mu_h$ and $\mu_\ell$, and the variance functions by $v_h(\cdot)$ and $v_\ell(\cdot)$, respectively; also $\mu := \mu_h + \mu_\ell$ (where $\mu < c$) and $v(\cdot) := v_h(\cdot) + v_\ell(\cdot)$. We note that in this priority setting we cannot restrict ourselves to centered processes. We denote the amount of traffic from the $i$th hp source in $(s,t]$, with $s < t$, by $A_{h,i}(s,t)$; we define $A_{\ell,n}(s,t)$ analogously. Also $\Gamma_h(s,t), \Gamma_\ell(s,t)$ and $R_h, R_\ell$ are defined as before.

REMARK 5.1. Notice that this setting also covers the case that the number of sources of both classes are *not* equal. Assume, for instance, that there are $n\alpha$ lp sources. Multiplying $\mu_\ell$ and $v_\ell(\cdot)$ by $\alpha$ and applying the fact that the Normal distribution is infinitely divisible, we arrive at $n$ i.i.d. sources.

In the tandem situation we could, without loss of generality, center the Gaussian sources. It can be checked easily that such a reduction property does not hold in the priority setting, since there is no counterpart of Remark 2.6. Hence we cannot assume without loss of generality that $\mu_h = \mu_\ell = 0$.

Analogously to Lemma 2.4, we obtain that $\mathbb{P}(Q_{\ell,n} > nb)$ equals

$$\mathbb{P}\left(\exists\, t>0 : \forall\, s>0 : \frac{1}{n}\sum_{i=1}^{n} A_{h,i}(-t,-s) + \frac{1}{n}\sum_{i=1}^{n} A_{\ell,i}(-t,0) > b + c(t-s)\right).$$

Let $J_p$ be the exponential decay rate of $\mathbb{P}(Q_{\ell,n} > nb)$; analogously to Theorem 3.1 it can be shown that this decay rate exists. Similarly to the tandem case, with $f(\cdot) \equiv (f_h(\cdot), f_\ell(\cdot))$,

$$\mathcal{S}_p^{s,t} := \{f \in \Omega \times \Omega : f_h(-s) - f_h(-t) - f_\ell(-t) > b + c(t-s)\},$$

$$\mathcal{U}_p^{s,t} := \left\{ f \in \Omega \times \Omega : \begin{array}{l} -f_h(-t) - f_\ell(-t) \geq b + ct; \\ f_h(-s) - f_h(-t) - f_\ell(-t) \geq b + c(t-s) \end{array} \right\}, \tag{26}$$

$$\mathbb{P}(Q_{\ell,n} > nb) = \mathbb{P}\left(\left(\frac{1}{n}\sum_{i=1}^{n} A_{h,i}(\cdot); \frac{1}{n}\sum_{i=1}^{n} A_{\ell,i}(\cdot)\right) \in \bigcup_{t>0}\bigcap_{s>0} \mathcal{S}_p^{s,t}\right).$$

THEOREM 5.2. *The following lower bound applies:*

$$J_p \geq \inf_{t>0} \sup_{s>0} \inf_{f \in \mathcal{U}_p^{s,t}} I(f), \tag{27}$$

*with* $\bar{f}_h(t) := f_h(t) - \mu_h t$, $\bar{f}_\ell(t) := f_\ell(t) - \mu_\ell t$ *and*

$$I(f) := \tfrac{1}{2}\|\bar{f}_h\|_{R_h}^2 + \tfrac{1}{2}\|\bar{f}_\ell\|_{R_\ell}^2.$$

The infimum over $f \in \mathcal{U}_p^{s,t}$ can be computed explicitly, as in Lemma 3.4. As the analysis is analogous to the tandem case, but the expressions are more complicated, we only sketch the procedure. Again there is a regime in which one of the two constraints is redundant. Define

$$k_p(s,t) := \mathbb{E}(A_h(s) | A_h(t) + A_\ell(t) = b + ct).$$

Using the convexity of the large-deviations rate function, it can be shown that, if

$$\mathbb{E}(A_h(t-s) + A_\ell(t) | A_h(t) + A_\ell(t) = b + ct) > b + c(t-s),$$

only the first constraint in (26) is tightly met; it is equivalent to require that $k_p(s,t) < cs$. [If $k_p(s,t) \geq cs$, either both constraints in (26) are met with equality, or only the second constraint is met with equality; exact conditions for these two cases are easy to derive, but these are not relevant in this discussion.] As before, under $k_p(s,t) < cs$, we obtain the decay rate

$$\inf_{f\in\mathcal{U}_p^{s,t}} I(f) = \frac{(b+(c-\mu)t)^2}{2v(t)}, \tag{28}$$

compare the FIFO queue with link rate $nc$; in the other cases the expressions are somewhat more involved. Denote by $t^F$ the value of $t>0$ that minimizes the right-hand side of (28).

Similarly to the tandem case, there is a regime (i.e., a set of values of the link rate $c$) in which $J_p$ coincides with the decay rate of an FIFO queue. In this regime, which we call regime (A), conditional on a large value of the total queue length, it is likely that the hp queue is empty, such that all traffic that is still in the system is in the lp queue. Hence, for all $c$ in

$$\{c|\forall\, s>0 : k_p(s,t^F) < cs\} \tag{29}$$

we conclude

$$J_p = \frac{(b+(c-\mu)t^F)^2}{2v(t^F)}.$$

If $c$ is not in the set (29), we can use the methodology of Section 3 to find a condition under which the lower bound of Theorem 5.2 is tight; we call this regime (B).

Remark 5.3. In the tandem case, we found that the FIFO result holds for $c_1 \geq c_1^F$, whereas it does not hold for $c_1 < c_1^F$; the threshold value $c_1^F$ was found explicitly in Section 3.3. In the priority setting there is not such a clear dichotomy. Consider, for instance, the situation in which both types of sources correspond to Brownian motions; $v_h(t) \equiv \lambda_h t$, $v_\ell(t) \equiv \lambda_\ell t$ and $\lambda := \lambda_h + \lambda_\ell$. Define

$$\Xi := \sqrt{\mu_\ell^2 + \frac{\lambda_\ell}{\lambda_h}(c-\mu_h)^2}.$$

Then straightforward calculus yields that for $(\lambda_h - \lambda_\ell)c \leq \lambda_h(\mu_h + 2\mu_\ell) - \lambda_\ell\mu_h$, regime (A) applies (i.e., the FIFO result holds):

$$J_p = \frac{2b(c-\mu)}{\lambda},$$

whereas otherwise we are in regime (B):

$$J_p = \frac{b(\Xi - \mu_\ell)}{\lambda_\ell};$$

this is shown by verifying that the lower bound of Theorem 5.2 is tight for the specific case of Brownian motion input. Using $\mu_h + \mu_\ell < c$, it can be verified easily that this implies that for $\lambda_h \leq \lambda_\ell$ the FIFO solution applies, whereas for $\lambda_h > \lambda_\ell$ only for

$$c \leq \frac{\lambda_h(\mu_h + 2\mu_\ell) - \lambda_\ell \mu_h}{\lambda_h - \lambda_\ell},$$

the FIFO solution applies.

5.2. *Discussion.* Large deviations for priority queues have been studied in several papers. We mention here the work by Mannersalo and Norros [23] and Wischik [33]. We briefly review their results, and compare them with our analysis. Our lower bound then reads

$$J_p^{(\mathrm{I})} := \inf_{t>0} \sup_{s>0} \Upsilon_p(s,t) \qquad \text{with } \Upsilon_p(s,t) := \inf_{f \in \mathcal{U}_p^{s,t}} I(f).$$

Just as we did, Mannersalo and Norros [23] identify two cases. They get the same solution for our regime (A), that is, the situation in which, given a long total queue length, the hp queue is relatively short; see also Berger and Whitt's [6] empty buffer approximation.

In regime (B) the hp queue tends to be large, given that the total queue is long. To prevent this from happening, [23] proposes a heuristic that minimizes $I(f)$ over

$$\{f \in \Omega \times \Omega : \exists\, t > 0 : -f_h(-t) - f_\ell(-t) \geq b + ct; -f_h(-t) \leq ct\}. \tag{30}$$

Because regime (B) applies, the optimum paths in the set (30) are such that the constraints on $f$ are tightly met; consequently (30) is a subset of $\mathcal{U}_p^{t,t}$. Hence the resulting decay rate, which we denote by $J_p^{(\mathrm{II})}$, yields a lower bound, but our lower bound will be closer to the real decay rate:

$$J_p^{(\mathrm{II})} := \inf_{t>0} \Upsilon_p(t,t) \leq \inf_{t>0} \sup_{s>0} \Upsilon_p(s,t) = J_p^{(\mathrm{I})}.$$

REMARK 5.4. In the simulation experiments performed in [23], the lower bound $J_p^{(\mathrm{II})}(b)$ is usually close to the exact value. Our numerical experiments (cf. the examples on the tandem queue in Section 4) show that the hp buffer usually starts to fill shortly after the total queue starts its busy period. This means that in many cases the error made by taking $s = t$ is relatively small. It explains why the heuristic based on set (30) performs well.

Wischik [33] focuses on discrete time, and allows more general traffic than just Gaussian sources. Translated into continuous time, in regime (B), his lower bound on the decay rate $J_p^{(\mathrm{III})}$ (Theorem 14) minimizes $I(f)$ over

$$\begin{aligned} \{f \in \Omega \times \Omega : \exists\, t > 0 : \exists\, s > 0 : -f_h(-t) - f_\ell(-t) \geq b + ct;& \\ -f_h(-s) \leq cs\};& \end{aligned} \tag{31}$$

again a straightforward comparison gives that our lower bound $J_p^{(\mathrm{I})}$ is closer to the actual decay rate:

$$J_p^{(\mathrm{III})} := \inf_{t>0} \inf_{s>0} \Upsilon_p(s,t) \le \inf_{t>0} \sup_{s>0} \Upsilon_p(s,t) = J_p^{(\mathrm{I})}.$$

REMARK 5.5. Recent work by Mannersalo and Norros [24] suggests that a similar approach could work for a queue operating under the Generalized Processor Sharing (GPS) scheduling discipline. For two classes of traffic (both with $n$ sources), sharing a resource with link capacity $nc$ and two buffers, the model is parametrized by the *weights* $\phi_1, \phi_2 \in [0,1]$, summing to 1. If both queues are nonempty, both classes receive their *guaranteed service rates* $n\phi_1 c$ and $n\phi_2 c$, respectively. If one class does not use all its bandwidth, it can be taken over by the other class in a work-conserving manner. For more details on the system mechanics for GPS, see, for instance, [28, 29].

Consider the probability that the first queue exceeds level $nb$, under the assumption that the mean input rates of both classes are smaller than their respective guaranteed service rates. Notice that the backlog of type 2 does not exceed that of an FIFO queue with link rate $n\phi_2 c$. This suggests that, in self-evident notation, the decay rate is well approximated by the infimum of $I(f)$ over

$$\{f \in \Omega \times \Omega : \exists\, t > 0 : \forall\, s > 0 : -f^{(1)}(-t) - f^{(2)}(-t) - f^{(2)}(-s) \\ \ge b + ct - c\phi_2 s\}.$$

Reasoning heuristically (see also [24]), it is not likely that (i) queue 2 is nonempty at the start of the busy period preceding overflow of queue 1, (ii) there is traffic left in queue 2 at the epoch queue 1 reaches overflow. This would lead to a minimization over

$$\{f \in \Omega \times \Omega : \exists\, t > 0 : \forall\, s > 0 : -f^{(1)}(-t) - f^{(2)}(-t) \ge b + ct; \\ -f^{(2)}(-s) \le c\phi_2 s\},$$

compare the sets $\mathcal{U}^{s,t}$ (as identified for the tandem system) and $\mathcal{U}_p^{s,t}$ (priority system). A lower bound for this decay rate is again found by taking the infimum over $t > 0$ and the supremum over $s > 0$, as before.

## APPENDIX

PROOF OF THEOREM 3.1. To prove Theorem 3.1, we first present an auxiliary result (cf. [27], Proposition 4.2). The set $\mathcal{S}$ is open; we now determine its closure. Define

$$\mathcal{S}^t := \{f \in \Omega : \forall\, s \in (0,t) : f(-s) - f(-t) > b + c_2 t - c_1 s\},$$

$$\mathcal{S}^{s,t} := \{f \in \Omega : f(-s) - f(-t) > b + c_2 t - c_1 s\}.$$

LEMMA A.1. *The closures of $\mathcal{S}^t$ and $\mathcal{S}$ are characterized as follows:*

$$\overline{\mathcal{S}^t} = \{f \in \Omega : \forall s \in (0,t) : f(-s) - f(-t) \geq b + c_2 t - c_1 s\}, \tag{32}$$

$$\overline{\mathcal{S}} = \bigcup_{t \geq t_0} \overline{\mathcal{S}^t}. \tag{33}$$

PROOF. We first prove (32). "$\subseteq$" is obvious:

$$\overline{\mathcal{S}^t} = \overline{\bigcap_{s\in(0,t)} \mathcal{S}^{s,t}} \subseteq \bigcap_{s\in(0,t)} \overline{\mathcal{S}^{s,t}}.$$

Now consider "$\supseteq$." Let $f$ be in the right-hand side of (32). Define, with $y^+ := \max\{0,y\}$ and $y^- := \min\{0,y\}$,

$$f_n(u) := f(u) + \frac{1}{n}(u^- + t)^+.$$

It is easy to see that (i) $\|f - f_n\|_\Omega \to 0$, and (ii) $f_n \in \mathcal{S}^t$; here (ii) follows from

$$f_n(-s) - f_n(-t) = f(-s) - f(-t) + \frac{1}{n}(t-s) > b + c_2 t - c_1 s$$

for $s \in (0,t)$. This proves (32).

Next we show (33). Again we establish two inclusions. "$\supseteq$" is done by picking an arbitrary $f$ from the right-hand side:

(i) Suppose there is a $t > t_0$ such that $f \in \overline{\mathcal{S}^t}$; then we can reuse the above argument: take an $f$ from the right-hand side of (33), and show that there is a sequence $f_n$ in $\mathcal{S}$ such that $\|f - f_n\|_\Omega \to 0$. This is exactly as before.

(ii) Suppose $f$ is only in the union in the right-hand side of (33) for $t = t_0$; then we have to show that $f$ can be approximated by an $f_n \in \mathcal{S}^{t_0+\delta_n}$, with $\|f - f_n\|_\Omega \to 0$, and $\delta_n := 1/n$. This is done by the following sequence:

$$f_n(t) = \begin{cases} f(t), & \text{for } t > -t_0, \\ c_1(t+t_0) + f(-t_0), & \text{for } t \in [-t_0 - \delta_n, -t_0], \\ -c_1\delta_n + f(t) + f(-t_0) - f(-t_0 - \delta_n), & \text{for } t < -t_0 - \delta_n. \end{cases}$$

Now $f_n \in \mathcal{S}^{t_0+\delta_n}$, as can be seen as follows. For $s \in (0, t_0)$, using $f \in \overline{\mathcal{S}^{t_0}}$ in conjunction with (32),

$$\begin{aligned} f_n(-s) - f_n(-t_0 - \delta_n) &= f(-s) + c_1\delta_n - f(-t_0) \\ &\geq b + c_2 t_0 - c_1 s + c_1 \delta_n > b + c_2(t_0 + \delta_n) - c_1 s, \end{aligned}$$

due to $c_1 > c_2$. For $s \in [t_0, t_0 + \delta_n)$, similarly,

$$\begin{aligned} f_n(-s) - f_n(-t_0 - \delta_n) &= -c_1 s + c_1 t_0 + c_1 \delta_n \\ &= -c_1 s + b + c_2 t_0 + c_1 \delta_n > b + c_2(t_0 + \delta_n) - c_1 s. \end{aligned}$$

Now concentrate on "$\subseteq$"; take $f \in \overline{\mathcal{S}}$:

(i) Hence there is a sequence $f_n \in \mathcal{S}$ such that $\|f - f_n\|_\Omega \to 0$. Because $f_n \in \mathcal{S}$, there is a sequence of epochs $t_n$ (all of them strictly larger than $t_0$) such that $f_n \in \mathcal{S}^{t_n}$.

(ii) $t_n$ is bounded. This can be seen as follows. Clearly, due to $f_n \in \mathcal{S}^{t_n}$,

$$f_n(-s) - f_n(-t_n) > b + c_2 t_n - c_1 s$$

for all $s \in (0, t_n)$. Hence also $-f_n(-t_n) \geq b + c_2 t_n$ (let $s \downarrow 0$), and consequently

$$\|f - f_n\|_\Omega \geq \frac{f(-t_n) - f_n(-t_n)}{1 + t_n} \geq \frac{b + c_2 t_n}{1 + t_n} + \frac{f(-t_n)}{1 + t_n}.$$

Suppose $t_n$ were not bounded; letting $n \to \infty$ would lead to a contradiction: $0 \geq c_2$ [use that $f(u)/(1+u) \to 0$ for $u \to \infty$].

(iii) Hence we can pick a subsequence $t_{n_k}$ such that $t_{n_k}$ goes to some finite limit $t_\infty \geq t_0$ for $k \to \infty$. Now $f \in \overline{\mathcal{S}^{t_\infty}}$, since, for all $s \in (0, t_\infty)$ and $k$ sufficiently large,

$$f_{n_k}(-s) - f_{n_k}(-t_{n_k}) \geq b + c_2 t_{n_k} - c_1 s.$$

This proves the lemma. $\square$

We now prove the theorem. Clearly, from Schilder's result,

$$\inf_{f \in \overline{\mathcal{S}}} I(f) \leq J \leq \inf_{f \in \mathcal{S}} I(f).$$

To show the stated, we prove that the infima over $\mathcal{S}$ and $\overline{\mathcal{S}}$ coincide.

Let $\mathcal{T}$ be defined by (12), and let $\mathcal{T}^t$ and $\mathcal{T}^{s,t}$ be defined analogously to $\mathcal{S}^t$ and $\mathcal{S}^{s,t}$; their closures are determined as in Lemma A.1. It is evident that the infima over $\mathcal{S}$ and $\mathcal{T}$ coincide. As mentioned above, our aim is to prove that the infima over $\mathcal{S}$ and $\overline{\mathcal{S}}$ match, but it turns out to be more convenient to show that the infima over $\mathcal{T}$ and $\overline{\mathcal{T}}$ match.

This is done by choosing $f$ from $\overline{\mathcal{T}} \cap R$ arbitrarily, and showing that we can approximate it by a path in $\mathcal{T}$. Clearly $f \in \overline{\mathcal{T}^{t^\star}}$ for some $t^\star \geq t_0$. Let $\zeta$ be an arbitrary path in $R$ that is strictly positive in $(0, t^\star]$, and define $f_n := f + \zeta/n$. Then there is a $t_n > t^\star$ such that $f_n \in \mathcal{T}^{t_n}$. This can be seen as follows.

(i) First observe that $f(s) \geq b + c_2 t^\star - c_1(t^\star - s)$ for all $s$ in the *closed* interval $[0, t^\star]$. Hence for all $s \in (0, t^\star]$, and $n \in \mathbb{N}$,

$$f_n(s) = f(s) + \frac{1}{n}\zeta(s) > b + c_2 t^\star - c_1(t^\star - s).$$

(ii) As this inequality also holds for $s = t^\star$, we conclude that there is a $t_n > t^\star$ with $f_n(s) > b + c_2 t_n - c_1(t_n - s)$ for all $s \in (0, t_n)$, or, equivalently, that $f_n \in \mathcal{T}^{t_n}$.

Now notice that, for $n \to \infty$,

$$\|f_n\|_R^2 = \left\| f + \frac{1}{n}\zeta \right\|_R^2 \to \|f\|_R^2,$$

which proves Theorem 3.1. □

Proof of Lemma 3.9. The lemma is proven in three steps. Notice that, as we are in regime (B), it holds that $c_1 < c_1^F = K(t_{c_2}^F)$, with

$$K(t) := \sup_{s\in(0,t)} \frac{k(s,t)}{s}. \tag{34}$$

(i) In [9], Lemma 3.1, it is shown that, under (5) and Assumption 3.3, $L_{c_2}(t)$ is decreasing for $t < t_{c_2}^F$, and increasing for $t > t_{c_2}^F$.

(ii) We now prove by contradiction that $K(t^\star) \geq c_1$. Suppose $K(t^\star) < c_1$. Then, by (34), for all $s \in (0,t^\star)$ it holds that $k(s,t^\star) < c_1 s$, and hence also $\sup_{s\in(0,t^\star)} L(s|t^\star) = 0$; see (17). Now consider the decomposition of Corollary 3.6:

$$\sup_{s\in(0,t)} \Upsilon(s,t) = L_{c_2}(t) + \sup_{s\in(0,t)} L(s|t); \tag{35}$$

$t^\star$ is minimizer of this expression. Because $k(s,t)$ is continuous, also for the *closed* interval $[0,t^\star]$ it holds that $k(s,t^\star) < c_1 s$. Hence it is possible to decrease $t^\star$ such that the first term in the right-hand side of (35) decreases (as we approach $t_{c_2}^F$ from above, see Step 1), while the second remains 0. Hence the sum of both terms decreases, implying that $t^\star$ cannot be optimal. So it cannot be that both

$$K(t^\star) < c_1 \quad \text{and} \quad t^\star > t_{c_2}^F.$$

Similarly $K(t^\star) < c_1$ rules out $t^\star < t_{c_2}^F$. Hence $K(t^\star) < c_1$ implies $t^\star = t_{c_2}^F$. However, $K(t_{c_2}^F) > c_1$. Contradiction.

(iii) Notice that $K(t^\star) \geq c_1$, in conjunction with (17), directly implies that $k(s^\star, t^\star) \geq c_1 s^\star$. This proves the lemma. □

Proof of Lemma 3.12. First we show that (23) holds. As noted earlier, both $m(\cdot)$ and $\rho(\cdot)$ have a maximum 1 at $-s^\star$. This means that (23) is necessary to enforce $m(r) \leq \rho(r)$ for $r$ in a neighborhood of $-s^\star$.

Next we show that (23) is equivalent to (24). First multiply both $m(\cdot)$ and $\rho(\cdot)$ by $h(\cdot)$, where

$$h(r) := \sqrt{\frac{\bar{v}(r)}{\bar{v}(-s^\star)}}(\bar{\mathbb{E}}A(-s^\star, 0) - c_1 s^\star).$$

Since $h(r)\colon(-t^\star,0)\to\mathbb{R}_+$, this yields the requirement $\pi(r)\le n(r)$ for all $r\in(-t^\star,0)$, with

$$\pi(r):=\bar{\mathbb{E}}A(r,0)+c_1r \quad \text{and} \quad n(r):=\frac{\bar\Gamma(r,-s^\star)}{\bar v(-s^\star)}(\bar{\mathbb{E}}A(-s^\star,0)-c_1s^\star).$$

Recall that $m(\cdot)$ and $\rho(\cdot)$ have the same function value and derivative at $-s^\star$. It is easy to derive that this implies that $(m\cdot h)(-s^\star)=(\rho\cdot h)(-s^\star)$ and $(m\cdot h)'(-s^\star)=(\rho\cdot h)'(-s^\star)$. Therefore, the necessary condition becomes $\pi''(-s^\star)\le n''(-s^\star)$.

Using standard formulas for conditional means of multivariate Normal random variables,

$$\bar{\mathbb{E}}A(r,0)=+\frac{\Gamma(-r,t^\star)}{v(t^\star)}(b+c_2t^\star),$$

leading to

$$\frac{d^2}{dr^2}(\bar{\mathbb{E}}A(r,0)+c_1r)\bigg|_{r=-s^\star}=\frac{b+c_2t^\star}{2v(t^\star)}(v''(s^\star)-v''(t^\star-s^\star)).$$

Assuming $r\le -s^\star$,

$$\bar\Gamma(r,-s^\star)=\frac{v(-r)+v(s^\star)-v(-r-s^\star)}{2}-\frac{\Gamma(-r,t^\star)\Gamma(s^\star,t^\star)}{v(t^\star)},$$

such that

$$\frac{d^2}{dr^2}\bar\Gamma(r,-s^\star)\bigg|_{r=-s^\star}=\frac{v''(s^\star)-v''(0)}{2}-\frac{v''(s^\star)-v''(t^\star-s^\star)}{2}\frac{v(s^\star,t^\star)}{v(t^\star)}.$$

It can be checked that the same result holds when the derivative is calculated for $r>-s^\star$. Now it is a straightforward but tedious computation to prove that this implies that $\pi''(-s^\star)\le n''(-s^\star)$ is equivalent to (24). □

**Acknowledgment.** The authors thank Krzysztof Dębicki (CWI) for valuable discussions.

CWI
P.O. Box 94079
1090 GB Amsterdam
and
Faculty of Electrical
Engineering, Mathematics
and Computer Science
University of Twente
P.O. Box 217
7500 AE Enshede
The Netherlands
e-mail: michel@cwi.nl

CWI
P.O. Box 94079
1090 GB Amsterdam
and
Faculty of Economics
and Business Administration
Vrije Universiteit
De Boelelaan 1105
1081 HV Amsterdam
The Netherlands
e-mail: miranda@cwi.nl